\documentclass[11pt,reqno]{amsart}
\usepackage{hyperref} 

\usepackage{that46}

\usepackage{graphicx,epsfig} % support the \includegraphics command and options
\DeclareGraphicsExtensions{.eps}
\title{}

\address{M.M.Graev, Scientific Research Institute for System Studies of Russian Academy of Sciences, Nakhimovsky prosp., 36, korpus 1, Moscow, 117218,  Russia}
\email{mmgraev@gmail.com}

\begin{document}

\def\contentsname{Contents}
\def\refname{References}
\def\bibname{References}
\def\figurename{Fig.}



%{\tiny\par\noindent \today \hfill THAT46.TEX 16:52 04.06.2012 (18.04.2012 (alkid-e.tex 1.04.2012))}

\begin{center}{}\Large\bf
On the compactness  of  the  set of  invariant Einstein metrics
\end{center}

\medskip

\begin{center}{}\large
M. M. Graev
\end{center}

\footnotetext{Supported by RFBR, grant 10-01-00041a.}

\begin{comment}
\medskip

\begin{center} 
Scientific Research Institute for System\\ Studies of Russian Academy of Sciences,\\ Nakhimovsky prosp., 36, korpus 1,\\ Moscow, 117218,  Russia
\end{center}

\begin{center} 
%\email{mmgraev@niisi.msk.ru}
%mmgraev@niisi.msk.ru, mmgraev@gmail.com
\texttt{mmgraev@gmail.com}
\end{center}
\end{comment}
%\null
\maketitle

\vskip-1cm
%\vskip-1cm
%\vskip-1cm
%\vskip-1cm
%\vskip-1cm



% A B S T R A C T

% \begin{quote}
{\small 
\noindent
\textsc{Abstract}.
\let\bar\overline
Let $M = G/H$ be a connected simply connected homogeneous manifold of a
compact, not necessarily connected Lie group $G$. We will assume that the
isotropy $H$-module $\mathfrak {g/h}$ has a simple spectrum, i.e. irreducible
submodules are mutually non-equivalent.

  There exists a convex Newton polytope $N=N(G,H)$, which was used for the
estimation of the number of isolated complex solutions of the algebraic
Einstein equation for invariant metrics on $G/H$ (up to scaling). Using the
moment map, we identify the space $\mathcal{M}_1$ of invariant Riemannian
metrics of volume 1 on $G/H$ with the interior of this polytope $N$.

  We associate with a point ${x \in \partial N}$ of the boundary a homogeneous
Riemannian space (in general, only local) and we extend the Einstein equation
to $\bar{\mathcal{M}_1}= N$. As an application of the Aleksevsky--Kimel'fel'd
theorem, we prove that all solutions of the Einstein equation associated with
points of the boundary are locally Euclidean.

  We describe explicitly the set $T\subset \partial N$ of solutions at the
boundary together with its natural triangulation.

  Investigating the compactification $\bar{\mathcal{M}_1} $ of $\mathcal{M}_1$,
we get an algebraic proof of the deep result by B\"ohm, Wang and Ziller about
the compactness of the set $ \mathcal{E}_1 \subset \mathcal{M}_1$ of Einstein
metrics. The original proof by B\"ohm, Wang and Ziller was based on a different
approach and did not use the simplicity of the spectrum.
 In Appendix we consider the non-symmetric K\"ahler homogeneous spaces $G/H$ with the second Betti number $b_2=1$. We write the normalized volumes $2,6,20,82,344$ of the corresponding Newton polytopes and discuss the number of complex solutions of the algebraic Einstein equation and the finiteness problem.

}
% \end{quote}

%   C O N T E N T S
\begin{quote}
\tableofcontents
\end{quote}
\vskip-1cm
\vskip-1cm

\clearpage

\begin{comment}
\begin{enumerate}{}%\let\se ction\item
\item{Invariant metrics on a compact homogeneous space $G/H$}
\item{Moment map and moment polytope}
\item{Compactification  $\Delta=\MET_1 \cup \Gamma $}
\item{Euclidean geometries at infinity} %{sect:4}
\item{Minimal compactification  $\Delta _{\min}$} %{sect:T}
\item{First application}
\item{Second application}
\item{Newton polytope and proof of Theorem~\ref{THM:e-nu}}
\item{Appendix. Case of K\"ahler homogeneous space with $b_2=1$}
\item{References}\dotfill\pageref{bibliography}
\end{enumerate}
\end{comment}

  

\section*{Introduction}

Let $M = G/H$ be  a  connected simply  connected homogeneous
manifold of a  compact, not necessarily connected  Lie group $G$. We will assume that the isotropy $H$-module $\mathfrak {g/h}$ has a simple spectrum, i.e. irreducible submodules  are mutually non-equivalent.

There exists a %compact 
convex Newton polytope $N=N(G,H)$, which was used for  the  estimation of  the number of isolated complex solutions of the algebraic Einstein equation for invariant metrics on $G/H$ (up to scaling), see \cite{2006,2007}.
Using the moment map, we identify the space $\mathcal{M}_1$ of invariant Riemannian metrics of volume 1 on $G/H$ with the interior of this polytope  $N$.

We associate  with a point ${x \in \Gamma = \partial N}$ of the  boundary  a  homogeneous Rieman\-nian space (in general, only local, since the stability subgroup can be non-closed)  and  we  extend the Einstein equation to  $\overline{\mathcal{M}_1}= N$. 
As an application of the Aleksevsky--Kimel'fel'd theorem,
we prove that all solutions  of  the Einstein  equation   associated with points of the  boundary are locally Euclidean.

We describe explicitly  the set $T\subset \Gamma$ of  solutions at the boundary together with its  natural  triangulation.
It is the standard geometric realization of a subcomplex of the simplicial complex, whose simplicies are the $H$-invariant subalgebras $\mathfrak t\subset\mathfrak g$ satisfying 
%$[\mathfrak t,\mathfrak t] \subset \mathfrak h \subsetneq \mathfrak t$ 
$\mathfrak t=\mathfrak h \oplus\mathfrak a$, $\mathfrak a\ne0$, $[\mathfrak a,\mathfrak a]=0$
(quasi toral subalgebras),  and the vertices are minimal such subalgebras.

Investigating the  compactification $\overline{\mathcal{M}_1} $ of $\mathcal{M}_1$, we get an algebraic proof of the deep result by B\"ohm, Wang and Ziller about the compactness  of  the  set  $\mathcal{E}_1 \subset \mathcal{M}_1$ of Einstein metrics.
 The  original proof by B\"ohm, Wang and Ziller \cite{BWZ} was based on a different approach and did not use the  simplicity of the spectrum.

In \S\ref{sect:1}, \ref{sect:2}, \ref{sect:3} 
% we give the definitions of 
we define the moment map and the moment polytope $\Delta$ (more general than $N$), and construct the corresponding compactification of $\MET_1$. Lather in \S 3 we prove that all solution of the Einstein equation at the boundary of $\Delta$ are Ricci-flat and, consequently, flat. 
In \S\ref{sect:4} we describe a triangulation of the set $T\subset\partial\Delta$ of these solutions. 
% We exhibit 
We consider examples, where $T$ is a finite set,
or a disjoint union of simplicies, or the join of two finite sets (a complete bipartite graph).
In \S\ref{sect:T} we construct the minimal (under inclusion) moment polytope $\Delta_{\min}$, 
% under inclusion of polytopes, 
by passing, if necessary, to some 'non-essential' extension of the coset space $G/H$. We prove that $T=\varnothing$, if the groups $G$ and $H$ are connected, and $\Delta=\Delta_{\min}$ (Proposition~\ref{PROP:C}).  
In \S\ref{sect:6} we use $\Delta_{\min}$ to prove the compactness of $\mathcal E_1$. We deduce it from the compactness of $T\cup\mathcal E_1$. 
It follows that $\mathcal E_1$ is compact 
in the case when the groups $G$ and $H$ are connected.
We sketch a proof in the general case.  
% In \S\ref{sect:7} we consider an application of the polytope $\Delta_{\min}$ to the non-existence problem for non-isolated complex solutions of the Einstein equation. We write the optimal upper bound of the normalized volume of $\Delta_{\min}$ for the number $\varepsilon$ of isolated solutions.  
% In \S\ref{sect:8} we outline that $\Delta_{\min} = N$, and get an upper bound for $\varepsilon$ (of the normalized volume of some permutohedron $\Pi \supset \Delta$).

In \S\ref{sect:7} we consider an application of the polytope $\Delta_{\min}$ to the finiteness problem for complex solutions of the algebrac Einstein equation.
We write the optimal upper bound of the normalized volume $\nu$ of $\Delta_{\min}$ for the number $\varepsilon$ of \textit{isolated} solutions.  
In \S\ref{sect:8} we outline that $\Delta_{\min} = N$, and get an upper bound for $\varepsilon$ (of the normalized volume of some permutohedron $\Pi \supset \Delta$, which is a central Delannoy number $D\in \{3,13,63,321,1683,\dots\}$, cf. \cite{Stanley2}).

In Appendix we consider the non-symmetric K\"ahler homogeneous spaces $G/H$
with the second Betti number $b_2=1$. In this case 
$
2^{-1}\nu \in\{1,3,10,41,172\},
$
and we have $\varepsilon=\nu$ for $2^{-1}\nu \in\{1,3,10\}$.
By the recent calculations of I.Chrysikos and Y.Sakane \cite{SACOS} it implies that for ${G/H=E_8/T^1\cdot A_3\cdot A_4}$ all complex solutions 
% of the algebraic Einstein equation 
are isolated, and $\varepsilon = 81$, so that $\nu-\varepsilon = 82-81=1$. 
The missing solution with multiplicity $1$ 'escape to infinity'. We indicate the missing solution explicitly. 
% We explain this difference by 'escape to infinity' of a solution with multiplicity $1$, that we indicate explicitly. 
We discuss a reduction of the finiteness problem for complex solutions
% non-existence problem of non-isolated complex solutions 
in the case of 
$
G/H = {E_8/T^1\cdot A_4\cdot A_2\cdot A_1}
$
(based on calculation of some 'marked' faces of $\Delta$ and
consideration of a toric variety $\Delta^{\CC}$),  
and prove that $\varepsilon<\nu$, where $\nu = 344$.



\section{Invariant metrics on a compact homogeneous space $G/H$}\label{sect:1}

% unless otherwise stated.

Let $ G / H $ be a connected simply connected $ n $-dimensional homogeneous space of a compact Lie group $ G $, \  $\rho : H \to \mathrm{GL}(\frak{g/h})$
the isotropy representation with a finite kernel.

Let us denote by $\MET = \MET(G, H)$ the cone
of invariant Riemannian metrics $g$ on $G/H$ (or, equivalently, $\rho(H)$-invariant  Euclidean  scalar products in   $\mathfrak{g/h}$),
and by $\MET_1 = \MET_1(G, H)$
the hypersurface of metrics $g$ with volume $\opn{vol}_g(G/H)=1$.

We shall assume, unless otherwise stated, that the representation $ \rho $ has a simple spectrum, 
%i.e. is split into $ d $ pairwise inequivalent irreducible summands
i.e., $\rho$ decomposes as a direct sum of $ d \le n $ pairwise inequivalent irreducible representations
(e.g., as it is  in the case
$ \opn {rank} (G) = \opn {rank} (H) $).
Therefore,
$$
\MET(G,H)  = (\RR_{>0})^{d}.
$$
We suppose ${d>1}$, and fix an
 ${H\cdot Z_G (H^0) }$-invariant Euclidean scalar product
$g_1$ on $\frak{g/h}$,  $ g_1 \in \MET_1 $.



\section{Moment map and moment polytope}\label{sect:2}

We will define  the {\bf moment map}
$
\displaystyle
\mu:
\MET(G,H) \to\RR^{n-1} = \left(\mbox{\begin{tabular}{c}
% the space  of 
diagonal matrices of\\ order $ n $ with trace $ 1 $
\end{tabular}}\right)
$
  as  the gradient of the logarithm of a ``suitable'' positive homogeneous
function on the cone $ \MET (G, H) $.
It is not unique and in particular depends on an ($H$-invariant) reductive decomposition
$$
\mathfrak{g}= \mathfrak{h} + \mathfrak{m}.
$$
For an invariant definition,
we  may   chose  the $B$-orthogonal decomposition, where $B$ is  the Killing  form (with the kernel $\mathfrak{z(g)}$).

We define a suitable function on $\MET(G,H) $ as the following modified scalar curvature
$$
\ell_{\theta }(G/H,g)= \opn{trace}(-(1+ \theta )\mathit{Ric}_{G/H,g}-B_{G/H,g}),
$$
where $\mathit{Ric}_{G/H,g}$ is the Ricci operator of a metric $ g $ at the point $eH \in G/H$, \,
$B_{G/H,g}= g^{-1} B |_{\frak m} \in \mathrm{End}(\mathfrak{m})$,
and $\theta$ is a parameter, $|\theta|<1$.

The corresponding moment map is given by
$$
\mu_{\theta }(g) = \frac{1}{\ell_{\theta }(G/H,g)}%^{-1}
(-(1+ \theta )\mathit{Ric}_{G/H,g}-B_{G/H,g}).
$$
Clearly, $\mu_{\theta}(g)$ belongs to
$\RR^{d-1}$, the space of $H$-invariant diagonal matrices with trace $1$.

\medskip

Changing $\mathfrak{m}$  to any other
$H$-invariant  
% $H \cdot Z_G(H^0)$-invariant 
complement $\mathfrak{m}'$ to $\mathfrak{h}$  we  get  another ``suitable'' function and  another moment map
$
\boxed
{ \mu = \mu_{\theta}
: \MET\to \RR^{d-1}, }
% : \MET_1\to \RR^{d-1}, }
$
which we call compatible with $\frak m$.
Now  we fix   any such  complement $\mathfrak{m}$  (not necessary  $B$-orthogonal).

\medskip

{\bf Remark} There are other possibilities to define a ``suitable'' function, but the scalar curvature
$$
\SC(G/H,g)=\opn{trace}(\mathit{Ric}_{G/H,g} )
$$
is not always a ``suitable'' function, since it can take non-positive values.
The following statements hold for the  'moment map'  $\mu$ associated with a suitable function  (more general than $\mu_{\theta}$).

\medskip

Now we associate with $\frak m$ a compact convex %polytope
polyhedron  $\Delta \subset \RR^{d-1}$.
Let
$$
\mathfrak{m}= \mathfrak{m}_1+\dots+\mathfrak{m}_d,
$$
where $\mathfrak{m}_i$ are irreducuble $H$-submodules of $\mathfrak{m}$.
Let $\varepsilon _i$, $i=1, \dots ,d$, be the weight
of the Lie algebra $\RR^d \subset \mathfrak{gl(m)}$
of $H$-invariant diagonal matrices,
such that $AX=\langle \varepsilon _i,A \rangle\,X$ for all $A \in \RR^d$, $X \in \frak m_i$.
By $\Delta $ we denote the convex hull of all weights of the form
\begin{center}{}
$
\varepsilon _i + \varepsilon _j - \varepsilon _k,
$
and $\varepsilon _r$,
\end{center}
where $g_1([\frak m_i,\frak m_j]+ \frak h,\frak m_k+ \frak h) \ne \{0 \}$, $B(\frak m_r,\frak m_r) \ne \{0 \}$
(we assume here $g_1(\frak h,\frak g):=0$).

\begin{EXAMPLES}{}\label{EXAM:21}
Let $G/H$ be a direct product of $d\ge 2$ isotropy irreducible spaces,
e.g., copies of $\CP^1$. Then $\Delta $ is the standard ${(d-1)}$-dimensional
coordinate simplex with vertices
$\varepsilon _1 = {(1, \dots ,0)}, \dots ,\varepsilon _d = {(0, \dots ,1)}$.
\end{EXAMPLES}

\begin{EXAMP}{}\label{EXAM:22}
Let $G/H$ be $SU(3)/T^2$.
Then $d=3$, and $\Delta $ is the triangle with vertices $(1,1,-1)$, $(1,-1,1)$, $(-1,1,1)$.
This valid for the spaces $G/H$ with $d=3$ and $[\frak m_i,\frak m_i] \subset \frak h$,
$[\frak m_i,\frak m_j] = \frak m_k$, $\{i,j,k \} = \{1,2,3 \}$. 
%%%assume $[\frak h,\frak m]= \frak m$ for later reference.
\end{EXAMP}

\begin{EXAMP}{}\label{EXAM:23}
Let $G/H$ be $E_8/(A_2)^4$ or $E_7/T^1\cdot (A_2)^3$. Then $d=4$,
and $\Delta $ is a $3$-polytope  with eight $2$-faces. It is an Archimedean solid
(a truncated tetrahedron), or respectively
the convex hull of two opposite faces of such a solid (a hexagon and a triangle).
\end{EXAMP}

Using a technical lemma from \cite[\S4.2]{Fulton}, one can prove the following theorem:

\begin{THM}\label{THM:DEL}

Let $G/H$ be a connected simply connected homogeneous space of a compact Lie group $G$ such that
$\frak{g/h}$ is a multiplicity-free $H$-module,
$ \MET_1 = \MET_1 (G, H) $ the space of the invariant Riemannian metrics of volume $ 1 $,
\   $\frak m$ an invariant complement to $\frak h$ in $\frak g$, 
\   $\mu : \MET_1\subset\MET\to \RR^{d-1}$ a moment map, 
compatible with $\frak m$, 
and $ \Delta \subset \RR ^{d -1} $ the compact convex polyhedron
associated with $\frak m$.

Then the map $\mu$ determines a diffeomorphism of the space $ \MET_1 $ onto the interior of $\Delta$.
We have $\dim (\Delta )=d-1$.

\end{THM}

We will consider the Euclidean %polytope 
polyhedron 
$\Delta\subset\RR^{d-1}$ as a compactification of the space $\MET_1$ of metrics, and call $\Delta$ the {\bf moment polytope} (associated with $\frak m$). The points on the boundary $ \Gamma = \partial \Delta $ of $\Delta$ we call {\bf points at infinity}.



\section{Compactification  $\Delta=\MET_1 \cup \Gamma $}\label{sect:3}

In this section, we associate with a point $x$ of the
boundary $\Gamma = \partial \Delta$ a Lie  algebra
 $$
 \mathfrak{g}_x = \mathfrak{h}+ \mathfrak{m}
 $$
  with a reductive decomposition and
  a \textit{fixed} $\mathrm{ad}({\mathfrak{h}})$-invariant Euclidean metric $g_1$ on $\mathfrak{m}$.
  Since, in general, the subalgebra $\mathfrak{h}$ generates a non closed
   subgroup  of the Lie  group $G_x$ associated  with $\mathfrak{g}_x$ , it does not define a  homogeneous Riemannian manifold.
   However,   we  can  exponentiate $\mathfrak{m}$  to a locally defined (non
  complete) Riemannian manifold $M^n(x)$
  with a transitive action of the Lie  algebra $\mathfrak{g}_x$  and  the stability subalgebra $\mathfrak{h}$.
  More precisely, we can speak about a germ of ``local  homogeneous Riemannain
  geometry''.
Later in this section we will describe explicitly the points at infinity corresponding to
germs of Einstein geometries.

\smallskip

To  describe the construction more carefully, we consider the
  moment map  $\mu : \MET_1\xrightarrow{\,\sim\,} \Delta\smallsetminus\Gamma$, and
 associate with each interior point $x=\mu(g) \in \Delta \smallsetminus \Gamma  $
 a homogeneous Riemanian space $(G/H,cg)$, $c>0$, where
  $cg \in \MET$ is a %unique
  Riemannian metric on $G/H$ (proportional to $g=\mu^{-1}(x)$) with
the same modified scalar curvature as $g_1$, namely,
$\ell_0(G/H,cg)=\ell_0(G/H,g_1)$. Here $ g_1 \in \MET_1 $ is the fixed $H\cdot Z_G (H^0) $-invariant Euclidean scalar product on
$\frak g/\frak h \cong \frak m$. 
Let $\phi : \frak m  \to G/H$ be a local diffeomorphism
 defined  in a neighbourhood  $M^{n}(x)$  of the origin by
$
\phi(Y)=\beta(a^{-1}Y)=\exp(a^{-1}Y)H
$
for all vectors  $ Y \in \frak m $ of sufficiently small length,
where
$a\in \mathrm{GL}(\frak m)$ is an
$ H $-invariant diagonal linear transformation on $ \frak m $ such that
$$
\ell_0(G/H,g)g(Y,Y) \equiv  \ell_0(G/H,g_1)g_1(aY,aY).
$$
We will consider $M^n(x) \subset \frak m $ as a Riemannian space
with respect to the metric $\phi^*(cg)$. Hence $\phi^*(cg)(Y,Y) =
g_1(Y,Y)$ for all tangent vectors in the origin $0\in M^n(x)$. We
  define  a   transitive Lie algebra of Killing vector fields on $M^n(x)$
$$
\frak g_x=(\frak h + \frak m, [\cdot,\cdot]_x)
$$
(isomorphic to $\frak g=Lie\, G=(\frak h + \frak m,[\cdot,\cdot])$)
 by $[T_aY,T_aZ]_x=T_a[Y,Z]$, where $T_aY=Y$ for all $Y\in\frak h$,
$aY$ for all $Y\in\frak m$. (So that
$[\cdot,\cdot]_{\mu(g_1)}=[\cdot,\cdot]$.) Clearly, the  stability
group  $H$ acts isometrically on $M^n(x)$.

We  will denote by $(M(x),g_1)$ and $M(x)$ the germs of the above
Riemannian homogeneous structure and,  respectively, the  homogeneous structure
on $M^n(x)$ at the point $0\in M^n(x)$. Non-formally,
$M(x)$  can  be  considered as a neighbourhood  $M^{n}(x)$
equipped with actions of $\frak g_x$ and $H$.

One can check that the Lie algebra $\frak g_x=(\frak h + \frak m, [\cdot,\cdot]_x)$
can be defined for every $x \in \Delta $ so that $\frak g_x$ depends continuously of $x$.
(However, $\frak g_x\not\cong \frak g$
for $x\in\Gamma$). In this way, the germ $(M(x),g_1)$ is
well-defined for all $x\in \Delta$.  Moreover, it satisfies the
following properties:
\begin{itemize}{}

\item
the Ricci tensor $\opn{ric}(M(x),g_1)$ and all others associated with the metric tensors
at the point $0\in M(x)$ depend continuously of $x$; 
cf. \cite{Lauret};

\item
the compact group $H$ acts on this germ isometrically with the fixed point $0$
and the same isotropy representation $\rho $ at $0$;

\item
the modified scalar curvature $\ell_0$ of a germ is well defined and
 is constant on $\Delta $, that is $\ell_0(M(x),g_1) =
\ell_0(G/H,g_1)$ for all $x \in \Delta $.

\end{itemize}

\begin{DEF*}{}

By {\bf infinitesimal homogeneous Riemannian space} we will understand
a quadruple
$(\mca{A}, \TIL{\frak g}, \TIL{\frak h}, g)$, where $\TIL{\frak g}$ is
a Lie algebra, $\TIL{\frak h}\subset \TIL{\frak g}  $ a subalgebra, $\mca{A}$
is a compact group of automorphisms of the pair $(\TIL{\frak g}, \TIL{\frak h})$,
 and $g$ is a Euclidean scalar product on  $\TIL{\frak g}/ \TIL{\frak h}$, invariant under
$\mca{A}$ and $\TIL{\frak h}$.

\end{DEF*}

Thus, we can associate with any point
$ x \in \Delta $
an infinitesimal homogeneous Riemannian space
$ (M (x), g_0) = (\mca{A},\frak g_x,\frak h, g_0) $, which we call also
a geometry, such that $\TIL{\frak h}
\cong \frak h$ and $\mca{A} \cong Ad_{\frak g}(H)$; the isomorphism
of Lie algebras $\TIL{\frak h}=\frak h$ extends to a isomorphism of
$\mca{A}$-modules
$\TIL{\frak g}$
and $\frak g = \frak h + \frak m$, which are identified,
and $g_0\in\MET$. 

%%%%%%%%% Finally, we can set $g_0=g_1$.

\medskip

A geometry $(M(x), g_1)$
associated with a point at infinity
$ x \in \Gamma $
%This  infinitesimal geometry
can  be exponentiated to a local
geometry $M^n(x)$, as above, but not necessary to a  global
homogeneous Riemannian geometry, since the  stability
subalgebra $\tilde{ \mathfrak{h}}=\mathfrak{h}$ can generate a non-closed  subgroup
(cf. Exam.~\ref{EXAM:53} below).

 However, we can apply to  such local homogeneous  Riemannian geometry
  the  Alekseevsky--Kimel'fel'd theorem, stating that the Ricci--flat homogeneous
  Riemannian geometries are locally Euclidean \cite {Al-Ki}.
Due to the fact that any Lie algebra $ \TIL {\frak g}= \frak g_x$ (which is a contraction of the
  compact Lie  algebra $\mathfrak{g}$) is
 of the type $(R)$, the  proof of the theorem  given in \cite {Al-Ki} can be modified
  so that it remains valid for a local homogeneous Riemannian manifold.
(The condition of simplicity of the spectrum of  $\rho$ is
insignificant.) Using this,
\boxed{we} prove
% one can prove
the following theorem:

\begin{comment}
\input THAT50/Theorem2
% \par\medskip
\sect ion{Euclidean geometries at infinity}\label{sect:4}
\input THAT50/Lemma1
\par\medskip
\end{comment}

\begin{THM}{}

Any  Einstein geometry at infinity   $(M(x), g_1)$ is locally Euclidean.

\end{THM}

\begin{proof}[Outline of proof] It is sufficient to prove
that the geometry $(M(x), g_1)$ is Ricci-flat, i.e., the scalar curvature $s=\SC(M(x), g_1)$ vanishes.

Let $\phi \subset \Gamma$ be any facet of the moment polytope $\Delta$ through the point $x$. Up to sign, there is a unique vector $z=(z_1,\dots,z_d) \in \ZZ^d$ with $\gcd(z_1,\dots,z_d)=1$, orthogonal to $\phi$, so $\langle x,z\rangle=0$.
We may assume that $z$ generates an edge of the following $d$-dimensional convex polyhedral cone\,:
$$
\nabla = \{y\in\RR^d : \langle x',y\rangle\ge0,\enskip \forall\,x'\in\Delta\},
$$
since otherwise we may pass from $z$ to $-z$.
Let $y,y',y''\in\RR^d$ and $y_i=\max(y'_i,y''_i)$
(respectively $\min(y'_i,y''_i)$) for all $i\in\{1,\dots,d\}$.
In this situation we write $y=\max(y',y'')$
and $y=\min(y',y'')$, respectively. We have
$$
\max(y,y')\in \nabla,\quad \min(y,0)\in \nabla,\quad 
\forall \,y,y'\in \nabla.
$$  
This follows from definitions of $\Delta$ and $\nabla$,
since $\bigoplus_{y_i<0} \frak m_i \subset \frak {z(g)}$.
Hence
$$
z=\max(z,0) + \min(z,0) = \max(z,0) \mbox{\ \ or\ }\min(z,0)
$$
(moreover, in the second case we have $\sum z_i = -1$).

We outline two proofs that $s=0$.
By Theorem~\ref{THM:DEL}, the moment polytope $\Delta=\ov{\mu_{\theta}(\MET_1)}$ is independent of $\theta$.
Using this, one can check that
for all $\theta\in (-1,1)$ 
$$
\phi\,\ni\,
\frac{1}{\ell_{\theta }(M(x),g_1)}
(-(1+ \theta )\mathit{Ric}_{M(x),g_1}-B_{M(x),g_1}).
$$
Therefore $z$ and ${r=\mathit{Ric}_{M(x),g_1}}$ 
can be considered as two orthogonal vectors in $\RR^d$,
%is orthogonal to $z$,
$\langle r,z\rangle=0$, 
and
$
s\,\sum z_i\dim \frak m_i = n \langle r,z\rangle=0,
$
so $s=0$.

For another proof of $s=0$, we may consider $z$ as a derivation of the Lie algebra $\frak g_x$ with the eigenspaces $\frak g_x^k$ such that $\frak h \subset\frak g_x^0$ and $\frak g_x^k\cap\frak m = \bigoplus_{z_i=k}\frak m_i$ (possibly, $\frak g_x^k=0$). Then
$$
\mbox{either\ \ }
\frak g_x=\bigoplus_{k=0}^{\infty} \frak g_x^k\,,
\mbox{\ \ or\ \ }\frak g_x=\bigoplus_{k=-\infty}^0 \frak g_x^k= \frak g_x^{-1} + \frak g_x^{0},
$$ 
and 
$[\frak g_x^k,\frak g_x^l]_x \subset \frak g_x^{k+l}$
for all integer $k,l$ (moreover, $\frak g_x^{-1} \subset \frak z(\frak g_x)$). Consider now $z$ as an element of $\frak {gl(m)}$, and assume $g^{\lambda} = e^{-\lambda z} .\,g_1$ is the one-parametric family of Euclidean scalar products on $\frak m$ (so that $g^\lambda\in \MET$ and $g^0=g_1$). 
We conclude that the geometries $(M(x), g^\lambda)$ with fixed $x$ and all $\lambda \in \RR$ are equivalent, and, hence, Einsteinian with the same scalar curvature $s$. 
By Hilbert--Jensen theorem \cite{Jen-2},  \ 
$$
0=\frac{d}{d\lambda} ((\det g^\lambda)^{1/n}\, s)
= s \,\frac{d}{d\lambda} e^{2\lambda\opn{trace}(z)/n}.
$$
(This is correct, since $\frak g_x$ is the Lie algebra of an unimodular Lie group. Note also that the Hilbert--Jensen theorem remains valid for a local homogeneous Riemannian manifold.) 
But $\opn{trace}(z)=\sum z_i \dim\frak m_i \ne0$, and, hence,
$s:=\SC(M(x),g_1)=0$.
\end{proof} 

\textbf{Remark.} The cone $\nabla$ is a ``tropical ring''
under operations $y\oplus y' =\max(y,y')$ and $y\odot y'=y+y'$, so that $y\odot(y'\oplus y'')=(y\odot y')\oplus(y\odot y'')$.



\section{Euclidean geometries at infinity}\label{sect:4}

Now we describe the points at infinity corresponding to locally Euclidean geometries.

\begin{LEM}{}\label{LEM:T}

Let $x \in \RR^d$. Then $x$ lies  in $\Gamma $ and the corresponding geometry $(M(x), g_1)$ is locally Euclidean if and only if $x$ belongs to  the convex hull of a subset of weights $\{\varepsilon _i
: i \in I \} \subset  \{\varepsilon _1, \dots , \varepsilon _d\}$,
such that the subspace $\frak m_I = \bigoplus_{i \in I} \frak m_i \subset \frak g$ satisfies conditions
$$
[\frak m_I,\frak m_I]=[\frak m_I,\frak h]=\frak m_I \cap \frak {z(g)}=0.
$$

\end{LEM}

\begin{proof}[Proof] 
Assume $\mu=\mu_{\theta}$. 
Then 
$
\xi = \frac{1}{\ell_{\theta }(M(\xi),g_1)}
(-(1+ \theta )\mathit{Ric}_{M(\xi),g_1}-B_{M(\xi),g_1})
$
for all interior points $\xi$ of $\Delta$ and, hence,
for all poits $\xi$ of the boundary $\Gamma$. 

\smallskip

Suppose that $x\in\Gamma $ and the corresponding
geometry $(M(x), g_1)$ is locally Euclidean. Then
\begin{equation}
x=\frac {B_{M(x),g_1}}{\opn{trace}  B_{M(x),g_1}}  
= \sum_{i=1}^d t_i\varepsilon_i,
\tag{*}
\end{equation}
for some coefficients $t_i\ge0$ with $\sum t_i =1$. Let
$$
\frak m_{\tau}= \bigoplus _{t_i>0} \frak m_i,\quad
\frak n= \bigoplus _{t_i=0} \frak m_i
$$
and let $\xi$ be a relative interior point of the 
convex hull $\tau$ of the set $\{\varepsilon_i : t_i>0\}$,
e.g., the point $\xi = x$. It is easy to check that $\frak m_{\tau}\cap \frak {z(g)} =0$, and $\xi \in \Delta$. It follows from $x \in \Gamma$ that $\xi \in \Gamma$.

We prove that %$\xi \in \Gamma$ and 
the corresponding geometry $(M(\xi),g_1)$
also is locally Euclidean, assuming that $\frak m$ is a subalgebra
of $\frak g_x$, i.e., $[\frak m,\frak m]_x\subset\frak m$.
(For example, if the reductive decomposition $\frak g = \frak h + \frak m$ is $B$-orthogonal, and, hence, $B_{\frak g_x}(\frak h,\frak m)=0$, then undoubtedly $[\frak m,\frak m]_x\subset\frak m$, since the stability subalgebra $\frak h$ contains a maximal semisimple subalgebra of the Lie algebra $\frak g_x$.)

Note that a necessary and sufficient condition for  $\frak m = \frak m_\tau + \frak n$ to be a transitive effective Lie algebra of motions of the Euclidean space $(\frak m, g_1)$ is  
\begin{equation}
\begin{array}{cc}
{[\frak m_\tau,\frak m_\tau]'}=[\frak n,\frak n]'=0,
&
[\frak m_\tau,\frak n]' \subset \frak n,
\\[1ex] 
g_1([\frak m_\tau,Y]',Y) =0,
% = g_1(\frak m_\tau,Y) =0,
&
\forall\, Y \in\frak n.
\end{array}
\tag{**}
\end{equation}
where $[\cdot,\cdot]'$ is the commutator on $\frak m$, 
$[Y,Z]'=[Y,Z]_x$. (Cf., e.g.,  \cite[\S5]{Jen-2}.)
This is clear, since 
$\frak n$ is the kernel of the  Killing form of the Lie algebra $(\frak m,[\cdot,\cdot]')$, by construction, and
\begin{equation}
g_1(\frak m_\tau,\frak n) =0.
\tag{***}
\end{equation} 
% (Cf. Hano, Jensen \cite{Jen-2}, and Heintze descriptions of left-invariant Riemannian metric on a solvlable Lie group with Ricci curvature or scalar curvature zero.) 
Define now $a=(a_1,\dots,a_d)\in (\RR_{>0})^d$ by $\sum a_i^2t_i \varepsilon_i = \xi $ and $a_j=1$, if $t_j=0$. Consider $a$ as an $H$-invariant diagonal linear transformation of $\frak m$, so that $a|_{\frak m_i}=a_i$, $i=1,\dots,d$. Define a new commutator $[\cdot,\cdot]'$ on $\frak m = \frak m_\tau + \frak n$ with 
the property \thetag{**} 
% properties \thetag{**}, \thetag{***} 
by $[Y,Z]' = a^{-1}[aY,aZ]_x$. Let $\frak m^{a}=(\frak m,[\cdot,\cdot]')$ be the corresponding Lie algebra, and $\exp(\frak m^a)$ the associated simply-connected Lie group. 
So $\exp(\frak m^a)$ is a metabelian group, and the scalar product $g_1$ on $\frak m$ 
gives
% induces 
the left-invariant Euclidean metric on  $\exp(\frak m^a)$.
The Killing form of $\frak m^a$ is $$B'=(\opn{trace}(B_{M(x),g_1}))\xi = (\ell_0(G/H,g_1))\xi .$$ Then there $\exp(\frak m^a)$ is locally equivalent to $(M(\xi),g_1)$, so that $[Y,Z]'=[Y,Z]_\xi$ for all $Y,Z\in\frak m$, and the assertion follows.

\smallskip

We prove now that $[\frak m,\frak m]_x\subset\frak m$,
and $[\frak m_\tau,\frak m_\tau]=[\frak m_\tau,\frak h]=0$.
There exist two $H$-invariant diagonal matrices $A,A_0 \in \frak {gl(m)}$ such that
\begin{itemize}
\item 
$x=\lim_{\lambda\,\to\,+\infty} \mu(e^{-\lambda A}e^{-A_0}.\,g_1)$ (where $a.g(X,X)=g(a^{-1}X,a^{-1}X)$),
\item
$[\cdot,\cdot]_x = \lim_{\lambda\,\to\,+\infty} T_{e^{A_0}e^{\lambda A}} [T_{e^{-\lambda A}e^{-A_0}}(\cdot),T_{e^{-\lambda A}e^{-A_0}}(\cdot)]$ 
\end{itemize}
(cf. \cite[\S2.3]{Fulton}). Then $\frak g_x = \frak h + \frak m^0 + A\frak m$, where $\frak m^0 = \{X\in\frak m:AX=0\}$. 
The compactness of the group $G$ implies that 
the subspace $A\frak m$ is a nilpotent ideal of the Lie algebra $\frak g_x=(\frak h + \frak m, [\cdot,\cdot]_x)$, its complement $\frak g_x^0 := \frak h + \frak m^0 $ is a subalgebra of $\frak g_x$, and, moreover, $\frak g_x^0$
decomposes as a direct sum of 
its center $\frak a = \frak z(\frak g_x^0)$ and a compact semisimple subalgebra $\frak k$.
Further, $\frak g_x$ is a transitive subalgebra of the complete Lie algebra $\frak{so}(n)+\RR^n$ of motions of Euclidean space, by assumption, and 
$\frak h = \frak g_x \cap \frak {so}(n)$.
So $\frak h$ contains a maximal semisimple subalgebra of $\frak g_x$. Therefore, $\frak k \subset\frak h$ and 
$\frak m^0 \subset \frak a$, so $[\frak h+\frak m^0,\frak m^0]_x=0$. This proves that $[\frak g_x,\frak m]_x \subset \frak m$. Let $X\in [\frak h+\frak m^0,\frak m^0]$.
Then $\lim_{\lambda\,\to\,+\infty} T_{e^{\lambda A}}X=0$.
Thus $X\in \frak{z(g)}$ and, hence, $X=0$. This proves that
$[\frak h+\frak m_\tau,\frak m_\tau] =0$.

\smallskip

Suppose now that $I\subset\{1,\dots,d\}$ and
$[\frak m_I,\frak m_I]=[\frak m_I,\frak h]=\frak m_I \cap \frak {z(g)}=0.$ Let $P\in\frak{gl(m)}$ be the orthoprojector with the kernel $\frak m_I := \bigoplus _{i\,\in\, I} \frak m_i$. Obviously, the Lie operation 
$[X,Y]'=\lim_{\lambda\,\to\,+\infty} T_{e^{\lambda P}}[T_{e^{-\lambda P}}X,T_{e^{-\lambda P}}Y]$ on $\frak h + \frak m$ is well-defined. 
% and the condition \thetag{**} holds for 
% Then the properties \thetag{**}, \thetag{***} are satisfied
Then the property \thetag{**} is satisfied
 for $[\cdot,\cdot]'$, 
$\frak m_{\tau}:=\frak m_I=\{X\in\frak m:PX=0\}$, and $\frak n:=P\frak m$
since the scalar product $g_1$ on $\frak g/\frak h \cong\frak m$ is $Z_G(H^0)$-invariant. 
There is a point $x\in\Gamma$ such that $x=\lim_{\lambda\,\to\,+\infty} \mu(e^{-\lambda P}.\,g_1)$.
% (where $a.g(X,X)=g(a^{-1}X,a^{-1}X)$). 
We have $[\cdot,\cdot]_x = T_c^{-1}[T_c\,\cdot,T_c\,\cdot]'$ for some scalar operator $c$ on $\frak m$. 
It follows from \thetag{**}, \thetag{***}
that the geometry $(M(x),g_1)$ is locally Euclidean. 
Then the point $x$ has the form \thetag{*} with $\{i:t_i>0\}=I$; e.g., $x=\varepsilon_j$, if $I=\{j\}$. 
Hence, $x$ belongs to the relative interior 
% Hence, $x$ is a relative interior point 
of the convex hull of the set $\{\varepsilon_i: i\in I\}$.
This completes the proof of Lemma~\ref{LEM:T}.  
\end{proof}

\par\medskip

Let us denote by $T\subset \Gamma$ the set of the points  at infinity corresponding to locally Euclidean geometries:
$
\boxed{
T:=\{ t \in \Gamma :  \opn{Riem} (M(t),g_1) = 0\}.
}
$

Here are examples with non-empty set $T\subset \Gamma$ of locally Euclidean geometries at infinity.

\smallskip

Notations.
Define conjugate linear transformations $A$ and $B$
of $\CC^p = \bigoplus_{l \in \ZZ_p } \CC e_l$
by $Ae_l = ce_{l+1}$, and $Be_l= c \omega ^{l-1}e_l$,
where $\omega ^p=1$, 
$c=i^{\,1-p^2}$.
For $p=2$ and $3$ we have
\begin{align*}&
A=\left\|\begin{array}{cr}0&i\\i&0\end{array}\right\| , &&
B=\left\|\begin{array}{cr}i&0\\0&-i\end{array}\right\| ; &&
%\intertext{and}&
A=\left\|\begin{array}{ccl}0&0&1\\1&0&0\\0&1&0\end{array}\right\|, &&
B=\left\|\begin{array}{ccl}1&0&0\\0&\omega&0\\0&0&\omega^2\end{array}\right\|,
\end{align*}
where $\omega^3=1$.
For  $p\not\equiv 0 (4)$ we have
$A^p = B^p = (-1)^{p+1}E$,  \  $ABA^{-1}B^{-1}= \omega E$,
and $A, B$ generate a finite subgroup of $SU(p)$,
%
% Here p\ne 0(4)
%
% p=4k+2; A= -i S,  B = -i D, det(S)=det(D)= -1
%       ( S  \sim D = diag(1,..,omega^{p-1}))
% det(A) = det (B) = -(-i)^p = -(-1)^{2k+1} = + 1
%
which we denote $J_p \subset SU(p) $ and call Jordan's group.
So $J_2$ is the group of the quaternionic units.
Further,
$\frak{su}(p)$
is the direct sum of abelian subspaces
%\begin{center}{}
$
\frak m_{(k,l)} =\frak m_{(-k,-l)} = \frak {su}(p)\,\cap\, ( \CC A^{k}B^l+ \CC A^{-k}B^{-l}).
$
%\end{center}

\smallskip

The complete bipartite graph $K_{r,s}$ is the graph with
$r+s$ vertices $a_1, \dots ,a_r$ and $b_1, \dots ,b_s$, and 
with one edge between each pair of vertices $a_i$ and $b_j$ (so $rs$ edges in all).

\begin{EXAMPLES}{}\label{EXAM:31}
Let $G=SU(p) \lX (\ZZ_p)^2$, $p$ be a prime, and $H=(\ZZ_p)^2$,
where $(\ZZ_p)^2=J_p/\ZZ_p$ is the group of automorphisms of $SU(p)$
generated by $Ad(A)$ and $Ad(B)$.
Let $(k,l) = Ad(A^kB^l)$ for all $k,l \in \ZZ_p$.
The $H$-modules
$\frak m_h = \frak m_{-h}$, $h \in H/(\pm1)$, $h\ne (0,0)$
are irreducible, and pairwise non-equivalent.
By regarding $\ZZ_p$ as a field, we have
\begin{center}{}
$[\frak m_h, \frak m_{h'}]=0$ iff the vectors $h, h' \in(\ZZ_p)^2$ are proportional.
\end{center}
Therefore
the set $T$ corresponding to the coset space $G/H$
is a disjoint union of simplices $\sigma _i$, $i=1, \dots ,p+1$,
with $\dim (\sigma_i ) = \max(0,(p-3)/2)$.
%
%(with %empty intersections $\sigma _i\cap \sigma _j = \varnothing$, $i<j$).
%

Similarly, if $G/H$ is the direct product of two such spaces with ${p \in \{2,3 \}}$,
then
$\dim(T){=}1$,
and $T$ is one of bipartite graphs $K_{3,3}$, $K_{3,4}$, or $K_{4,4}$.
\end{EXAMPLES}

\begin{EXAMP}{}\label{EXAM:32}
For an example with $\frak h\ne 0$, take
$G=SU(p+q+1)\lX ( J_p \times J_q  )$
and  $H= T^2 \times J_p \times J_q$, where $p,q \in\{2,3 \}$,
and $T^2$ is a torus.
Then $G/H$ has a simple spectrum of the isotropy rep\-re\-sen\-ta\-tion,
and $T$ is the complete bipartite graph $K_{p+1,q+1}$.
\end{EXAMP}

%PROWERITX: class  {\bf a  } {\bf\large\checkmark}\marginpar{{\bf\checkmark}}
% Prowereno ?

\medskip

We will now describe the vertices $v$ of the polytope $\Delta$ which belong to $T$.

\begin{LEM}\label{LEM:v}

Let $(M(v), g_1)$, $v\in T$, be a locally Euclidean geometry at infinity.
Assume for simplicity that the fixed decomposition $\frak g=\frak h +\frak m$ is $H\cdot Z_G(H^0)$-invariant.
Then the point $v$ is a vertex of the moment polytope $\Delta$, if and only if $v=\varepsilon_j$ for some $j$, and
$$
[\frak m_i,\frak m_j]\subset \frak m_i,
\qquad \forall \quad  i\in\{1,\dots,d\}.
$$

\end{LEM}

\begin{proof}[Proof]
Let $v$ be a vertex of $\Delta$. By Lemma~\ref{LEM:T}, \
$v=\varepsilon_j$, where
\begin{equation}
[\frak m_j,\frak m_j]=[\frak m_j,\frak h]=0,
\quad
[\frak m,\frak m_j]\ne 0.
\tag{*}
\end{equation}
Conversely, suppose a weight $\varepsilon_j$ satisfies (*).
Then $g_1([\frak m_i,\frak m_j],\frak m_k) = g_1([\frak m_k,\frak m_j],\frak m_i)$ for all $i,k$ 
% since $\frak m$ and $g_1$ are $Z_G(H^0)$-invariant.
since $g_1$ is $Z_G(H^0)$-invariant. 
Therefore $\varepsilon_j$ is either a unique vertex of $\Delta$ with $x_j>0$, or a half-sum of two distinct points
$p,q\in \Delta$,
$p\ne q$
of the form
\begin{align*}
&
p=\varepsilon_i+\varepsilon_j-\varepsilon_k,
&&
q=\varepsilon_k+\varepsilon_j-\varepsilon_i,
&&
i\ne j\ne k\ne i.
\end{align*}
In the second case, $\varepsilon_j$
is not a vertex, since $\Delta$ is convex.
In the first case, we obtain $g_1([\frak m_i,\frak m_j],\frak m_k)=0$ for all $k\ne i$
because $g_1([\frak m,\frak m_j],\frak m_j)=0$.
%%So $[\frak m_i,\frak m_j]=0$ or $\frak m_i$.
Lemma~\ref{LEM:v} follows.
\end{proof}



\section{Minimal compactification $\Delta _{\min}$}\label{sect:T}

Let $G/H$ be a connected simply connected homogeneous space of a compact Lie group $G$ such that
%$\frak{g/h}$ splits as a direct sum of $d>1$ irreducible submodules
$\frak{g/h}$ is a multiplicity-free $H$-module with at least two irreducible submodules,
%with ${d>1}$ irreducible submodules,
$ \MET_1 = \MET_1 (G, H) $ the space of the invariant Riemannian metrics of volume $ 1 $,
\   $\mu : \MET_1\to \RR^{d-1}$ the moment map,
%$\frak m\subset\frak g$ the corresponding complement of $\frak h$,
and $\Gamma$ the boundary of the %moment 
%$(d-1)$-dimensional
polytope $\Delta=\ov{\mu(\MET_1)}$.
%and $\Delta=\MET_1\cup\Gamma$ the corresponding compactification of $\MET_1$.
So $\dim \Gamma = d-2 \ge 0$.
The points 
%at infinity
$x\in\Gamma$ corresponds to geometries at infinity $(M(x), g_1)$.

\smallskip

The subset $ T \subset \Gamma $ of all locally Euclidean geometries at infinity 
(described in Lemma~\ref{LEM:T} above)
has a  natural triangulation, as the following  theorem states\,:

\begin{THM}\label{THM:T}

The  set  $ T \subset \Gamma $ of locally Euclidean geometries at infinity
is a union of some (closed) faces of the %standard coordinate 
$ (d-1) $-dimensional simplex $ S \subset \RR ^{d-1} $
with vertices $\varepsilon_i$, $i\in\{1,\dots,d\}$.

\end{THM}

In this section, we minimize this union $T$ by changing the moment map $\mu:\MET_1\to\RR^{d-1}$ and minimizing  
the moment polytope $\Delta=\ov{\mu(\MET_1)}$.
Moreover, we consider the maximal $T_{\max}$ and $\Delta_{\max}$
of $T$ and $\Delta$ (under inclusion).
Each $T$ is the union of all simlices of $T_{\max}$ that lie in $\Delta$.
The aim is to obtain the following compactification of $\MET_1 : $

\begin{DEF*}{}

A compactification $ \Delta=\MET_1\cup \Gamma $ of the space $\MET_1= \MET_1 (G, H) $
is called {\bf admissible} if $ T $ contains no whole faces of the boundary $ \Gamma $.

\end{DEF*}

In the case of an admissible compactification, one can check that $ \dim (T) <d-2 $ and, moreover, for each proper face
$ \gamma $ of the polytope $ \Delta $, we have
$$
\dim(T\cap \gamma ) < \dim (\gamma ).
$$

The map $ \mu $ and the moment polytope $ \Delta $ are defined with some freedom.
It depends on the  reductive decomposition $\mathfrak{g}= \frak h +  \frak m $.
There is a unique maximal moment polytope $\Delta_{\max}$,
containing all the others.
Its corresponds to  the $Q$-orthogonal  reductive decomposition,
where $ Q $ is any $Ad(G)$-invariant Euclidean metric on $ \frak g :$
$$
Q(\frak h, \frak m)=0.
$$
Although such complement $ \frak m $ looks %%to be
%
%%%%%%%%%%%%%%%%%%%%%%%%%%%%%%%%%%%%%%%%%%%%%%%%%%
\footnote{
moreover, $ \Delta _{\max} $ allows  to deal  only  with global  homogeneous geometries instead of local ones.
}
%%%%%%%%%%%%%%%%%%%%%%%%%%%%%%%%%%%%%%%%%%%%%%%%%%
most elegant and symmetric (cf. \cite {BWZ}),
it can give rise to a non-admissible compactification.
This holds, if and only if the set $T_{\max}$ of locally Euclidean geometries 
at the boundary of $\Delta_{\max}$ contains a vertex of $\Delta_{\max}$.

\medskip

Let us denote by $\Delta_{\min}$ the convex hull of all vertices $v$ of $\Delta_{\max}$ that do not lie in $T=T_{\max}$, and all vertices $v=\varepsilon_j$ of $S$ satisfying the same property $v\notin T_{\max}$.

\medskip

Turning to the spaces $G/H$ 
in the five examples above, we  have $\Delta _{\min} = \Delta _{\max}$,
but for the $(2k+1)$-dimensional sphere $U_{k+1}/U_{k}$, $k>0$, we have distinct segments
\begin{align*}
\Delta_{\max}&=[2\varepsilon_2-\varepsilon_1,\varepsilon_1],
&
\Delta_{\min}&=[2\varepsilon_2-\varepsilon_1,\varepsilon_2].
\end{align*}

\begin{OBS*}

A compactification $\Delta =  \MET_1\cup \Gamma $ is admissible iff $\Delta =\Delta_{\min}$.

\end{OBS*}

It is easy to check that $\Delta_{\min}$ is contained in all the moment polytopes $ \Delta $, but may be different from any of them
(e.g., for the sphere $G/H=SU_{k+1}/SU_{k}$, $k>1$).
If $\Delta_{\min}$ is a moment polytope,
its corresponds to  the $B$-orthogonal  reductive decomposition, that is,
$$
B(\frak h, \frak m)=0.
$$
Moreover, %the polytope $\Delta_{\min}$
its
depends only of the subspace
$\MET \subset \underline{\otimes^2T^*}(G/H)$
(cf. Proposition~\ref{PROP:Del_min=Nw}, below).

\medskip

We will show that extending the  group $G$ so that the space  $\MET_1(G,H)$ does not change, we can always construct an admissible compactification.
Suppose $G_1$ is a compact Lie group, the semidirect product of $G$ and a $G$-invariant torus:
\begin{center}{}
$G_1= (S^1)^k \opl\lX\limits_{\pi_0(G)} G$,
where
$(S^1)^k \subset \opn{Isom}(G/H, g_1)$,
and
${k\ge 1}$.
\end{center}
Assume, moreover, that $G_1$
acts almost effectively on the manifold $G/H$ (in a natural way)
with an isotropy subgroup $H_1 \supset H$.
So
\begin{center}{}
$G_1/H_1=G/H$, and $\dim(G_1)>\dim(G)$.
\end{center}
In this situation, we call the homogeneous space $G_1/H_1$ a {\bf toral extension} of the space $G/H$.
We call such extension {\bf non-essential}, if $\MET(G_1,H_1) = \MET(G,H)$,
and {\bf essential}, otherwise.

\begin{LEM}\label{LEM:T-2} The following conditions are equivalent:
\begin{enumerate}

\item
all toral extensions of $G/H$ are essential,
and $\frak{z(g)} \subset\frak m$, e.g.,
$
B(\frak h, \frak m)=0;
$
\item
 $T$ contains no vertices of $\Delta$, i.e.,
$\MET_1\cup \Gamma$ is an admissible compactification.
\end{enumerate}

\end{LEM}

We may assume that all toral extension of $G/H$ are essential, 
since one can always pass from $G/H$ to a (unique)
maximal non-essential toral extension of $G/H$,
which can be described explicitly.
Thus $\frak{z(g)} \subset\frak m$ iff
$\Delta$ contains in any other moment polytope, and
$$
\Delta = \Delta _{\min}.
$$
For example, this assumption is fulfilled, if
$\bigcap_{g \in \MET(G,H)} \opn{Isom} (G/H,g) = G $.

Remark that a toral extension of a connected group is also connected.

\begin{PROP}{}\label{PROP:C}

Suppose $G$, $H$ are connected groups, and $\Delta = \Delta _{\min}$.
Then there are no locally Euclidean geometries at infinity, that is, $T= \varnothing$.

\end{PROP}

\begin{proof}[Proof] By Lemma~\ref{LEM:T},  \  \    $T_{\max}$ is the empty set or a point.
By Lemma~\ref{LEM:v}, this point is a vertex of $\Delta_{\max}$.
Then $T=\varnothing$.
\end{proof}

Now we turn to examples with  $G$, $H$ connected,
where $\Delta _{\min} \ne \Delta _{\max}$.

\begin{EXAMPLES}{}

Consider the homogeneous space $M_{k,l}^{m,n} = (S^{2m+1} \times S^{2n+1})/T^1$
of $G={(U_{m+1} \times U_{n+1})/T^1}$
studied by Wang and Ziller (1990), see also \cite {B-G}.
The isotropy representation $\rho $ has a simple spectrum
($k,l,m,n>0$). %$kl\ne0$).  % CORRECTED !
%
% Pri $m=n=1, k+l=0  -- tozhe prostoj spektr, poskol'ku $G$ ne poluprosta.
%
Then $ \Delta _{\max} $ is a triangle, $ T_{\max} $ is one of its vertices, and $ \Delta _{\min} $ is a trapezoid, obtained by truncation of the triangle at the vertex $ T_{\max}: $
$$
\Delta _{\max} = \opn{Conv}\{(2,0,-1),(0,2,-1),(0,0,1)\}, \quad T_{\max}=\{(0,0,1)\},
$$
$$
\Delta _{\min} = \opn{Conv}\{(2,0,-1),(0,2,-1),(0,1,0),(1,0,0)\}.
$$
\end{EXAMPLES}
\begin{EXAMP}{}

Let $ M_{k, l} ^7 $ be a seven-dimensional homogeneous Aloff--Wallach space
with $k>l>0$
(so $\rho $ has a simple spectrum).
Then $ \Delta _{\min} $ is a (irregular)  octahedron. The polytope $ \Delta _{\max} $ has seven faces and seven vertices.
It can be obtained by constructing a tetrahedron on a face of $ \Delta _{\min} $.
The seventh vertex is $T_{\max}$.

\end{EXAMP}



\section{First application}\label{sect:6}

In \cite {BWZ}, the following theorem about the structure  of the set of invariant Einstein metrics on a compact homogeneous space was derived from a certain  variational theorem.

\begin{THM}{}\label{THM:1}

Let $G$ be compact Lie group, $G/H$ a
connected simply connected (or with finite fundamental group), homogeneous space, and $\mathcal E_1=\mathcal E_1(G,H)$ the set of all invariant, positive definite Einstein metrics on $G/H$
with
% unit
volume $ 1 $.
Then
$\mathcal E_1$ consists of at most finitely
many compact linearly connected components.

\end{THM}



The  set  $ \MET_1(G,H)$   of
all invariant unit volume Riemannian metrics $ g $  on $ G / H $,
$\opn{vol}_g(G/H) =1$,
has  the  structure of non-compact Riemannian
symmetric  space.
The  subset of Einstein metrics is
the set of critical points of an algebraic function,
 assigns to every metric $ g \in \MET_1 (G, H) $ the scalar curvature $ s = \SC (G / H, g) $, and, moreover, 
its gradient at $g$ is the minus traceless part of the Ricci tensor of $g$, that is,
for all $g\in\MET_1,$
$$ 
\opn{grad} s(g)=-\opn{ric}^0(g)
$$
 (Theorem of Hilbert--Jensen \cite {Jen-2,AB}).
 Therefore, Theorem~\ref{THM:1} is equivalent to the following proposition.

\begin{PROP}{}\label{PROP:1}

The subset $\mathcal E_1(G,H) \subset  \MET_1(G,H)$ is bounded.

\end{PROP}

As we shall see, the admissible compactification
$ \Delta _{\min} = \MET_1 \cup \Gamma $
leads to a simple, new, mostly algebraic, proof of these results for the special case of a homogeneous space with simple spectrum of the isotropy representation (i.e., in the case when all $ H $-invariant quadratic forms on $ \frak g / \frak h $ can be reduced  simultaneously  to principal axes).

\medskip

Remark that in the  original proof by B\"ohm, Wang and Ziller  \cite{BWZ} the  simplicity of the spectrum  was not used.



\begin{proof}[Outline of the proof for the case of a simple spectrum]
The set $ \mca {E} = T \cup \mca {E} _1 $ of all points of the polytope $ \Delta $ %(finite or at infinity),
(possibly at infinity),
corresponding to Einstein geometries, is compact
(we do not dwell on the proof).

This implies that in the case when the groups $G$,  $H$ 
% A consequence of this is that if the groups $G$ and $H$
are connected, then the set
$\mca{E}_1=\mca{E}_1(G,H)$ is compact.
 Indeed, we may assume that $\Delta = \Delta_{\min}$,
 and use Proposition~\ref{PROP:C}.
Then there are no locally Euclidean geometries at infinity, that is, ${T= \varnothing}$. The assertion follows.

In the general case it is sufficient to check
that $ T $ is open in $ \mca {E} $.
This is obviously true for $d=1+\dim \MET_1=2$, and we may assume $d>2$.
Every point $ t \in T $ lies in the closure of  a suitable submanifold of the form
$ \MET_1 (G_1, H_1) \subset \MET_1 (G, H) $.
Here $G_1/H_1$ is a toral extension of the space $G/H$,
where $G_1 \supset G$, $H_1 \supset H$, $G_1/H_1 = G/H$.
In this case the homogeneous manifold $ G_1/H_1 $ represents the same simply connected manifold as $ G / H $, and also has a simple spectrum of the isotropy representation. 
Assuming $\frak g=\frak h+\frak m$ is an $H\cdot Z_G(H^0)$-invariant decomposition (e.g.,$B$-orthogonal), then $\frak g_1=\frak h_1 + \frak m$ is an $H_1$-invariant reductive decomposition of the extended Lie algebra $\frak g_1$.
It follows from Theorem~\ref{THM:DEL} that the moment map $ \mu $ 
(comatible with $\frak m$) 
defines a diffeomorphism of $ \MET_1 (G_1, H_1) $ onto a linear submanifold of the interior of the polytope $ \Delta = \ov{\mu(\MET_1)} $, that is, its interesction with an affine plane.
The condition $ \Delta = \Delta_{\min} $ of Section~\ref{sect:T} implies that this submanifold is proper, i.e.,
$$
\dim \MET_1(G_1,H_1)<d-1=\dim \MET_1(G,H).
$$
We may assume by induction on $d$ that the proposition holds for $G_1/H_1$, and we remark that the submanifold $\MET_1(G_1,H_1)$ is invariant under 
the gradient Ricci flow $\dot g = -\opn{ric}^0(g)$ on $ \MET_1 (G, H) $.

\smallskip

Now we will associate with a point $t\in T$ an explicit submanifold $ \MET_1 (G_1, H_1)$ of the interior of $ \Delta $ described below. 
By Theorem~\ref{THM:T}, \  $T$ is the union of some faces of the standard weight simplex $S$ with vertices $\varepsilon_i$, $i\in \{1,\dots,d\}$. Let $\tau\ni t$ be the smallest face $\sigma$ of $S$ containing the point $t$, so that $\tau=\bigcap_{\,t\in\sigma\subset T}\sigma$, and let $\varepsilon_i$, $i \in I$, are vertices of $\tau$.
The corresponding submanifold $ \MET_1 (G_1, H_1) $ consists
of all interior points $x=(x_1,\dots,x_d)$ of the polytope
$\Delta$ satisfying the following system of linear equations:
$$
x_i=x_k, \quad \mbox{if\  }
g_1([\frak m_i,\frak m_j],\frak m_k)
=g_1([\frak m_k,\frak m_j],\frak m_i)\ne 0,
\mbox{\ for some\ }
j\in I.
$$
(Recall that $\sum x_i\dim \frak m_i = 1$ for all $x\in\Delta$.)

We can give an equivalent definition of $ \MET_1 (G_1, H_1)$. 
Denote by $\gamma=\bigcap \beta $ the intersections of all faces $\beta\subset\Delta$ such that $t\in \beta$ (so $t\in\tau\subset\gamma$).
Remark, that
$\sum_{k \notin I} x_k \dim \frak m_k =0$ for $x\in\gamma$, and $>0$ for $x\in\Delta\smallsetminus\gamma$. 
Moreover, 
since $\Delta=\Delta_{\min}$, this intersection $\gamma$ 
can be obtained explicitly as the convex hull of the points
\begin{align*}
&
\varepsilon_i+\varepsilon_j-\varepsilon_k,
&&
\varepsilon_k+\varepsilon_j-\varepsilon_i,
&&
j\in I
\intertext{where} 
&
{g_1([\frak m_i,\frak m_j],\frak m_k)\ne0},  
&&
i,k\notin I, 
&&
i\ne k.
\end{align*}
As we noted above, $\dim(\tau)<\dim(\gamma)$, since
$\Delta=\MET_1\cup \Gamma$
is an admissible compactification of $\MET_1$.
Consider $\RR^d$ as the Lie algebra of the group
$ (\RR_{> 0}) ^d \subset \mathrm{GL} (\frak {m}) $ with
the Euclidean metric $ (x, x) = \sum \dim(\frak m_i)x_i ^2 $ (so that $(\varepsilon_i,\varepsilon_j) = \frac{1}{\dim \frak m_i}\delta_{ij}$).
Let $\Omega$ be the sphere of unit vectors tangent to the face $\gamma$ and orthogonal to $\tau-t=\{z-t:z\in \tau\}$.
Let $Z$ be the intersection of $\Delta$ with
the orthogonal complement of the vector subspace $\opn{span}(\Omega)$ at the point $t$.
Then $Z$ is obviously a compact convex polytope 
of dimension $\ge 1$, containing the point $t$.
The intersection of $Z$ with the interior of $\Delta$
contains the point $\mu(g_1)$, and
 coincide\footnote{The compactification $Z$ of $ \MET_1 (G_1, H_1) $ is non-admissible, since $T$ contains the face $\tau$ of $Z$.}
 with $ \MET_1 (G_1, H_1) $, i.e.,
 $Z \smallsetminus \Gamma =  \MET_1 (G_1, H_1) $.

\smallskip

To carry out induction on $ d $, we must show that every Einstein metric $ g ^{G / H} \in \MET_1 (G, H) $, sufficiently close to $ t $ (if it exists) would be contained in $\MET_1(G_1,H_1)$.
We can define a small open neighborhood $U_{\rho}$ of the point $t$ in $\Delta$ by
$$
U_{\rho}=\{\lambda A+z: \lambda\in [0,\rho), A\in \Omega, z\in Z,|z-t|<\rho \}  .
$$
(By construction, it is an open subset of $\Delta$, if $0<\rho<\rho_0$).

\begin{LEM} The complement $U_{\rho} \smallsetminus Z$ contains no solution of the Einstein equation
%%%%%%%(possibly at infinity)
(that is, no point $ x\in\mca {E} = T \cup \mca {E} _1 $),
if $\rho$ is sufficiently small.
\end{LEM}

To prove this lemma, we consider the flat geometry $(M (t),g_1) $ as the geometry induced on a simply transitive group of motions of Euclidean space, and use the following facts.
The scalar curvature $ s (g) $ of each left-invariant Riemannian metric $ g $ on a solvable Lie group is non-positive, $ s (g) \le 0 $.
A metric $g_0$ with  $s (g_0) = 0 $  is Euclidean
(G.Jensen \cite {Jen-2}, E.Heintze),
and the Hessian $ s'' (g_0) $ of the function $g\mapsto s(g)$
has the rank $ = \opn {codim} \{g: s (g) = 0 \} $.
 We want to extend this Hessian  over  each geodesic line on $ \MET_1 (G, H) $
 orthogonal to $ \MET_1 (G_1, H_1) $
 (with respect to the natural inner Euclidean metric on $ \MET_1 (G, H) $).

More precisely, denote by $\SC(M(z),g)$ the scalar curvature of $(M(z),g)$, and consider $g\in \exp(\lambda \Omega).g_1$.
To each triple $ z \in Z $, $ A \in \Omega $, $ \lambda \ge0 $, we associate the number
$$
u(z,A,\lambda ) = -\frac12\frac{\partial }{\partial \lambda }\SC(M(z),\,e^{ - \lambda A}\opn{.}g_1).
$$
We have $\frac{\partial \textstyle u}{\partial \lambda }(t,A,0)>2 \delta >0$ for all $A \in \Omega $. 
This follows from the above facts about $s(g)=\SC(M (t), g)$, since $\frac{\partial \textstyle u}{\partial \lambda }(t,A,0)=-2s''(g_1)(g_1A,g_1A)$.
Moreover, $u(z,A,0) \equiv 0$ (in particular, when $ z \notin \Gamma $ this follows immediately from the invariance of $ \MET_1 (G_1, H_1) $
under the gradient Ricci flow $\dot g = -\opn{ric}^0(g)$ on $ \MET_1 (G, H) $).
Using continuity,  we get an estimation
 $
u(z,A,\lambda )\ge \delta \lambda , \,
\forall\, (z,A, \lambda ) \in Z' \times \Omega \times [0,\rho]
$
for a sufficiently small neighborhood $ Z'\subset Z $ of the point $ t $, and some  $ \rho> 0 $.

Now we can estimate the traceless part $ \opn {ric} ^0 $ of the Ricci tensor $ \opn {ric} $ at the point
$0\in M(x)$ for each of the infinitesimal Riemannian homogeneous spaces $(M(x),g_1)$  
with the parameter $x\in\Delta$ sufficiently close to $t$.
Changing $\rho$ if necessary, one can construct a natural locally one-to-one
continuous
map $\Phi: U_{\rho}\to \Delta$,
$
 (z, A, \lambda) \longmapsto z+\lambda A \in U_{\rho} \longmapsto x=x(z,A,\lambda ) \in \Delta ,
$
where $\lambda \in [0,\rho)$, possessing the following properties\,:
\begin{itemize}
\item
$
%\Phi(\lambda A+z) =
x(z,A,\lambda) =\mu(e^{-\lambda A}.\mu^{-1}(z)) %\qquad \forall
$
for all
$
(z,A,\lambda) \in U_{\rho}\smallsetminus \Gamma
$
(so the restriction $\Phi|_{U_{\rho}\smallsetminus \Gamma}$
can be considered as the normal
exponential map along $ \MET_1 (G_1, H_1) = Z\smallsetminus \Gamma $ with respect to the
$ (\RR_{> 0}) ^{d-1} $-invariant Euclidean metric on $ \MET_1 (G, H) =  \Delta \smallsetminus \Gamma$).

\item Moreover,
$\Phi|_{Z\cap U_{\rho}}=\opn{id} $, i.e.,
$x(z,A,0)\equiv z$.
For each face $\beta$ of $\Delta$ containing the point $t$
there is a smooth map 
$U_{\rho}\cap \mbox{relative interior}(\beta)\ni y \longmapsto \Phi(y)\in\mbox{relative interior}(\beta)$.
Every disc $D(z)=\{z+ \lambda A: A\in \Omega, \lambda \in[0,\rho)\} $, $z\in Z\cap U_{\rho}$ is tangent to $\Phi(D(z))$ at the center\,$z$.

\begin{comment} 
% \input THAT50/End_Proof-1 
\input THAT50/End_Proof-2
\end{proof}
\end{comment}

\item
$u(z,A, \lambda ) \ge \delta \lambda $,
for all points $z+\lambda A \in U_{\rho}$ (as above).
\end{itemize}
Consider now a scalar product $g\in\MET$, a point $x\in\Delta$, the Lie algebra $\frak g_x=(\frak h+\frak m, [\cdot,\cdot]_x)$,  and denote the geometry $(M(x),g)$  simply by $([\cdot,\cdot]_x, g)$. To any $H$-invariant linear transforma\-tion $a$ of $\frak m$ we associate a geometry $(a.[\cdot,\cdot]_x, a.g)$,
were $a.[\cdot,\cdot]_x = T_a[T_{a^{-1}}\cdot,T_{a^{-1}}\cdot]_x$ is a new Lie operation on $\frak h+\frak m$, and $a.g(\cdot,\cdot)=g(a^{-1}\cdot,a^{-1}\cdot)\in \MET$ is an $H$-invariant Euclidean scalar product on $\frak m$. 
% As an immediate consequence of the obvious equivalence between geometries $(a.[\cdot,\cdot]_x, a.g)$ and $([\cdot,\cdot]_x, g)$) we obtain the following Heber's identity (cf. \cite[\S6]{He}):
Geometries $(a.[\cdot,\cdot]_x, a.g)$ and $([\cdot,\cdot]_x, g)$) are equivalent, by construction.
As an immediate consequence we obtain the following Heber's identity (cf. \cite[\S6]{He}):
$$
\SC(a.[\cdot,\cdot]_x, a.g)=\SC([\cdot,\cdot]_x, g),
\qquad\forall \qquad a \in (\mathrm{GL}(\frak m))^H.
$$
% If $x=x(z,A,\lambda)$, then
For each  $x=x(z,A,\lambda)$
there is a scalar operator $\kappa$ on $\frak m$ 
such that $T_{\kappa^{-1}}[T_{\kappa}\cdot,T_{\kappa}\cdot]_x = e^{\lambda A} .\, [\cdot,\cdot]_z $. Then 
\begin{gather*} 
\langle \opn{ric}([\cdot,\cdot]_x, g_1),\, \kappa ^2A \rangle
= \left. - \frac12\frac{d}{dt}\right|_{\,t=0} 
\SC(e^{\lambda A} .\, [\cdot,\cdot]_z,\, e^{-tA} .\, g_1)
\\ 
= - \frac12\frac{\partial}{\partial\lambda}\SC([\cdot,\cdot]_z,\, e^{-\lambda A} .\, g_1)
= u(z,A, \lambda ). % \ge \delta \lambda .
\end{gather*}
Finally, for all $ x= x(z,A,\lambda)\in \Phi(U_{\rho}) $,
and some real function $ \kappa $ in $ x $ we have
$$
\langle \opn{ric}(M(x),g_1),\, \kappa ^2A \rangle
= u(z,A, \lambda ) \ge \delta \lambda .
$$
Clearly, $ \opn {trace} (A) = 0 $. This means that if $ \opn {ric} ^0 (M (x), g_1) = 0 $,
then $ \lambda = 0 $, and $ x \in Z $.
We have $U_{\rho'}\subset \Phi(U_{\rho})$ for some $\rho'>0$.
These imply Lemma and Proposition.
\end{proof}



A locally Euclidean geometry at infinity $ M (t) $
plays a central role in the above proof. 
Cf. the nice study of the flat space
$\RR^{n-k} \times T^k=\RR^{n-k} \times (\RR/\ZZ)^k$
as the limit of a sequence of compact homogeneous Riemannian spaces $(G_i/H_i, g_i)$ in \cite[\S2]{BWZ}.

\medskip

If \ $ T = \varnothing $, then the compactness of $ \mca {E}_1 $
is reduced to the compactness of $ \mca {E} $,
and the proof of the proposition is reduced to the first sentence.

\begin{EXAMPLES}{}
The condition $ T = \varnothing $ holds if $ \opn {rank} (G) = \opn {rank} (H) $
by Lemma~\ref{LEM:T}.
\end{EXAMPLES}

\begin{EXAMP}{}
%Let $G$ be a compact connected group, and $G/H$ a $S^1$-bundle over a K\"ahler homogeneous space $ G / K $, which is the subbundle of unit vectors of a untwisted ample line bundle. We call $G/H$ {\it a generalized Hopf bundle}.
%{\bf\large\checkmark}\marginpar{{\bf39}}

Let $G$ be a compact connected group, and $G/H$ the total space of a principal circle bundle over a K\"ahler homogeneous space $ G / K $, associated to a untwisted ample line bundle. We call $G/H$ {\it a generalized Hopf bundle}.

Let, moreover, the spectrum of the isotropy representation $\rho $ of
the group $H$ be simple.
Then the space $G/H$ has at most one toral extension (cf. Section~\ref{sect:T}).
If it exists, then $G$ is a semisimple group.
It follows from the simple spectrum condition, that
this extension $G_1/H_1$ is non-essential.
Passing from $G/H$ to $G_1/H_1$, we may assume that the space $G/H$ has no toral extension (so $\dim (\frak{z(g)})=1$).
Choose now the $B$-orthogomal complement $\frak m$ to $\frak h$, that is $B(\frak h, \frak m)=0$. Then $\Delta = \Delta _{\min}$. By Proposition~\ref{PROP:C}, \  $T=\varnothing$.
\end{EXAMP}

Note that the homogeneous spaces $ M_{k, l} ^{m, n} $ and $ M_{k, l}^7 $ considered
above are generalized Hopf bundles over ${\CP^{m}\times \CP^{n}}$ and ${F_3 (\CC) = SU_3 /T^2}$, respectively.   Assuming $B(\frak h, \frak m)=0$, then $ {T = \varnothing }$.

\medskip

Here is a simple example with $T=\varnothing$ and non-connected $G$, $H$.

\begin{EXAMP}\label{EXAM:53} %(S^3)^5
Let $G=(U_{k+1})^5 \lX C_5$, $H=(U_k)^5 \lX C_5$, so that
$G/H$ is the direct product of five spheres $S^{2k+1}=U_{k+1}/U_k$, and $C_5$ is the cyclic group of permutations of spheres. Then $d{=}4$. 
We have four irreducible $H$-modules $\frak m_i$, $i=1,\dots,4$ of dimensions $\dim \frak m_i = 1,2,2,10k$ respectively.
The polytope $\Delta_{\max}$ is an octahedron with vertices
$\varepsilon_i$, $\delta_i=2\varepsilon_4-\varepsilon_i$, $i=1,2,3$, and $\Delta_{\min}$ is a tetrahedron $(\delta_1, \delta_2, \delta_3, \varepsilon_4)$.  Let $B(\frak h,\frak m)=0$. Then $\Delta=\Delta_{\min}$. By Lemma~\ref{LEM:T}, \  $T=\varnothing$.
Moreover, an infinitesimal homogeneous Riemannian space $(M(x),g_1)$ 
is defined only locally (hence, is non-complete), if and only if $x$ is an 
%relative 
interior point of an edge $(\delta_r, \varepsilon_4)$, or a triangular
face $(\delta_1,\delta_r, \varepsilon_4)$,
where $r\in\{2,3\}$. 
In this case it is an isomorphism of Lie algebras $I_x:\frak g_x \cong \frak g$, but $I_x(\frak h)\ne \frak h$.
(Note that for the spaces $G/H$ in our other examples 
all geometries at infinity are locally isometric to complete Riemannian spaces.)
\end{EXAMP}



\section{Second application}\label{sect:7}

% NB  Besse conjecture with correction

 In \cite [Introduction] {BWZ} the authors asked the question about the  finiteness of the set $ \mca {E} _1 = \mca {E} _1 (G, H) $ of unit volume Einstein metrics  on a compact
 simply connected homogeneous space $ G / H $ with simple spectrum of isotropy representation.
 In this case, if $ d> 1 $, the Einstein equation reduces to a system of $ d-1 $
 rational algebraic equations on $ d-1 $ unknowns.
C.B\"ohm, M.Wang, and W.Ziller ask the following question: 
Is this system always {\it generic}, i.e.,
does it admit  at most finitely many {\it complex} solutions?

Here is a partial answer to this question \cite{2006,2007}.
Assume for the moment that $ \mca {E} _1 (\CC) $ is an infinite set.
Therefore it is noncompact, since the Einstein equation is algebraic.
Then it can be compactified by attaching some of the (complex) solutions at infinity lying on $ \Gamma^{\CC} $.
Here $ \Gamma^{\CC} $ can be regarded as a complex hypersurface in a compact complex algebraic variety with singularities $ \Delta ^{\CC} $, which is a complexification of the polytope $ \Delta $.
Thus, we obtain:

\begin{CLAIM*}

The set $ \mca {E} _1 (\CC) $ is finite if and only if in some neighborhood of the hypersurface $ \Gamma ^{\CC} $ all complex solutions are at infinity, i.e., lie on $ \Gamma ^{\CC} $.

\end{CLAIM*}

{\it Moreover, there are  no solutions at infinity,
if and only if $ \mca {E} _1 (\CC) $ is a finite set,
and, counting with multiplicities, it consists 
of
\begin{equation*}
\begin{aligned}{}
\nu
&
= (\opn{vol}(S))^{-1}\opn{vol}(\Delta_{\min})
\\ &
={(d{-}1)!}\opn{vol}(\Delta_{\min}).
\end{aligned}
\end{equation*}
solutions. } Here
$ \Delta = \Delta_{\min} $ is an admissible compactification of $ \MET_1 $,
and $S$ is the standard $(d-1)$-dimensional simplex in $\RR^{d-1}$.
(Note that $\nu$ is always an integer.)

\smallskip

These claims can be make rigorous and proved using the theory of toric varieties.
Since $\frak {g / h}$ is a multiplicity-free $H$-module,
the space of invariant Riemannian metrics on $ G / H $ has the natural complexification of the form
$$
\MET^{\CC} = (\CC \setminus 0)^d = (\CC \setminus 0)\times \dots \times (\CC \setminus 0).
$$
(Note that $ \MET ^{\CC} $ contains all the invariant pseudo-Riemannian metrics on $ G / H $.)
The quotient $\MET^{\CC}/\CC^\times $, where $\CC^{\times} = \{(z, \dots ,z)$, $z \in \CC \setminus 0\}$, can be considered as a complexification of the space $ \MET_1 $. 
The compactification $ \Delta = \MET_1 \cup \Gamma $ of the space $\MET_1 $ also has a natural complexification, namely, the toric variety
$$
\Delta^{\CC} = (\MET^{\CC}/\CC^\times) \cup \Gamma^{\CC}
$$
(see, e.g., \cite {Fulton}).
Here $\Delta^{\CC}$ is the toric variety of the fan in the lattice $\mathrm N$ from the polytope $\Delta=\Delta_{\min}\subset \varepsilon_d+\mathrm M_\RR$,
where $\mathrm M=\opn{Hom}(\mathrm N,\ZZ)=\sum_{i=1}^{d-1} \ZZ (\varepsilon_i-\varepsilon_{i+1})$.

The algebraic torus $(\CC\setminus 0)^d$ acts on $\Delta^{\CC}$ with open orbit $\MET^{\CC}/\CC^\times$, so that the subgroup $\CC^{\times}$ acts trivially.
In this way, 
 the polytope $\Delta$ can be considered as the closure of a single orbit of a subgroup $(\RR_{>0})^d$. The closure of each orbit of $(\CC\setminus 0)^d$ meets $\Delta$ in a single face, and each orbit of the compact torus $(S^1)^d$
meets $\Delta$ in a single point.

Algebraic Einstein equations are naturally defined on
$\MET^{\CC}$, $\MET^{\CC}/\CC^\times $, and $\Delta^{\CC}$. Let $\varepsilon = \varepsilon (G,H)$ be
the number of its {\it isolated} complex solutions (counting with multiplicities) on $ \MET^{\CC} / \CC^\times $.

Using the generalized Bezout theorem, one can prove the following theorem.

\begin{THM}{}\label{THM:e-nu}

Suppose $\frak {g / h}$ is a multiplicity-free module of $H$,
and $d=\dim \MET(G,H) > 1$.
Let $\nu = \nu(G,H) 
= (\opn{vol}(S))^{-1}\opn{vol}(\Delta_{\min})
={(d{-}1)!}\opn{vol}(\Delta_{\min})$.
Then
$$
\varepsilon \le \nu < 6^{\,d-1};
$$
for $ \nu = \varepsilon $ all solutions are isolated, and complex solutions at infinity cannot exist; 
for $ \nu> \varepsilon $ there is at least one complex solution lying on $ \Gamma ^{\CC} $ (i.e., at infinity).
\end{THM}

Roughly speaking, the missing $\nu-\varepsilon$ solutions ``escape to infinity''.

\smallskip

The strict inequality $ \nu> \varepsilon $ holds, for example,
if $ G $ is simple, $ H $ is its maximal torus, and
$G/H \ne SU(2)/T^1, \, SU(3)/T^2$
\cite{2006,2007}.

\smallskip

The estimation $ \nu \ge \varepsilon $ is sharp, as the following examples show.

\smallskip



\begin{EXAMPLES}{}
For $G/H = (S^{2k+1})^5$ in the preceding example,
we have $\nu=1$ and one obvious solution, so $ \nu = \varepsilon $.
For the spaces $G/H$
defined in
examples~\ref{EXAM:21}, \ref{EXAM:22} (if $[\frak h,\frak m]= \frak m$),
\ref{EXAM:23} ($G=E_7$ and~$E_8$)
hold $\varepsilon = \nu$.
This follows immediately %(without any calculations)
(without solving the Einstein equation)
from \cite[\S7.1, Tests 1 and 2]{2007},
Proposition~\ref{PROP:Del_min=Nw} below,
and the fact that $\Delta = \Delta _{\min}$.
By finding volumes, we obtain, respectively,
$$
\varepsilon = \nu = 1,1,4,20,23.
$$
Moreover, $SU(3)/T^2$ and every K\"ahler homogeneous space $G/H$  satisfying
conditions of Example~\ref{EXAM:22} admit $\nu=4$ positive definite invariant Einstein metrics 
$g\in\MET(G,H)$ with scalar curvature $1$
(D.V.Alekseevsky, 1987; M.Kimura, 1990).
\end{EXAMPLES}

\begin{EXAMP}{} For homogeneous spaces of Wang--Ziller and Aloff--Wallach, respectively, 
with $k,l,m,n>0$ and $k>l>0$
%(in the case of a simple spectrum of the isotropy representation $\rho $)
we have
$$
\varepsilon (M_{k,l}^{m,n}) = \nu (M_{k,l}^{m,n}) = 3,
\quad
\varepsilon (M_{k,l}^7) = \nu (M_{k,l}^7) = 16.
$$
To the proof,  one can examine faces  of polytopes
$ \Delta _{\min} $ using Tests 1 and 2 of \cite [\S7.1] {2007}.  Equality $ \varepsilon = \nu $ follows immediately without any calculations
%of the type of polyhedra $ \Delta _{\min} $ by Tests 1 and 2 of \cite [\S7.1] {2007}
(also in the second case, see \cite[Exam. 7.5]{2007}, is used the absence of complex Ricci-flat metrics on the underlying flag space  $SU(3)/T^2$).
\end{EXAMP}





\section{Newton polytope and proof of Theorem~\ref{THM:e-nu}}\label{sect:8}

In this section 
we interpret
$\Delta_{\min}$ as a Newton polytope,
we estimate the  normalized  volu\-me 
$\nu=(d-1)!\opn{vol} (\Delta_{\min})$,
and prove Theorem~\ref{THM:e-nu}.
Consider the moment polytopes $\Delta$ (and $\Delta_{\min}$) as polytopes
with vertices in $\ZZ^d$ by setting $\varepsilon_1=(1,0,\dots,0),\dots,\varepsilon_d=(0,\dots,0,1)$.

 We express a metric $ g \in \MET (G, H) $ as
$
g= \bigoplus x_k g_1 |_{\frak m_k},
$
and consider $x_k>0$, $ k \in \{1, \dots, d \}$,
as coordinates on $ \MET (G, H) $.
By ${s(g)=\SC(G/H,g)}$ we denote the scalar curvature of $g$.
Then 
$$
-\frac{x_i}{m_i}\frac{\partial s}{\partial x_i}
=
\frac{b_i}{2x_i}
- \frac{1}{4m_i} \opl{\textstyle\sum}_{j,k=1}^d {\scriptstyle [i,j,k]}
\frac{2x_k^2-x_i^2}{x_ix_jx_k},
\qquad 1\le i \le d,
$$
where $b_i>0$, ${\scriptstyle [i,j,k]}\ge0$ are coefficients, $m_i=\dim(\frak m_i)$, and
the original grouping of monomials taken from \cite {BWZ}.

The Einstein equation reduces to a system of $ d-1 $ homogeneous equations
$$
f_i:= \frac{x_i}{m_i}\frac{\partial s}{\partial x_i}-\frac{x_{i+1}}{m_{i+1}}\frac{\partial s}{\partial x_{i+1}} =0,
\quad 1\le i < d.
$$

\begin{figure}[ht]
\includegraphics{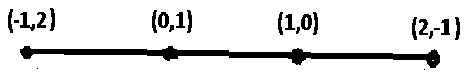}
% \includegraphics{3.eps}
%\includegraphics[width=65.0mm,height=13.8mm]{3.bmp}
%650 138
\caption{%Permutohedron $\Pi$. Case of $ d = 2$. Then $\Pi$
If $d=2$, then the permutohedron $\Pi$
 is a segment. It is equal to three segments $S=[(1,0),(0,1)]$,
and its normalized length is $P_1(3)=3$}
\end{figure}

\begin{figure}[ht]
\includegraphics{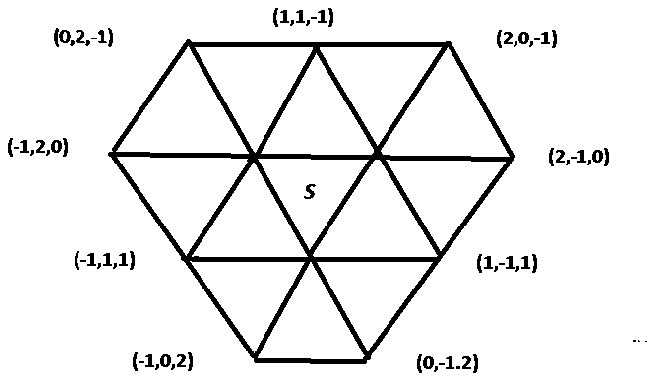}
% \includegraphics{13.eps}
%\includegraphics[width=68.3mm,height=45.6mm]{13.bmp}
%683 456
\caption{Case $ d = 3$. Then $\Pi$ is a hexagon. It is equal to thirteen triangles $S$,
and its normalized area is $P_2(3)=13$}
\end{figure}

We will use the following terminology. 
A Laurent polynomial in $x=(x_1,\dots,x_d)$ is a polynomial in $x_i^{\mathstrut},x_i^{-1}$, $i=1,\dots,d$. The {\bf Newton polytope} $\Nw(f)\subset\RR^d$ of a Laurent polynomial $f$ is the convex hull of the vector exponentials of its monomials. 

\smallskip

By considering $s$ and $f_i$ as homogeneous Laurent polynomials in $x^{-1}=(x_1^{-1},\dots,x_d^{-1})$, cf.\,\cite{2006}, we obtain:  

\begin{PROP}\label{PROP:Del_min=Nw} For a general linear combination $f=\sum c_i f_i$ of $f_i$ with coefficients $c_i\in \RR$ we have 
$
\Nw(f)=\Nw(s)=\Delta_{\min}.
$
\end{PROP}

Thus the polytope $ \Delta _{\min} $ can be introduced an invariant way as the Newton polytope $\Nw(s)$ of the polynomial of scalar curvature.

\begin{proof}[Outline of proof] Obviously $\Nw(f)\subset \Nw(s)$. Moreover, $\Nw(f)\ne\varnothing$. From the expression of $s$ in \cite[Eq. (7.39)]{AB} follows that 
$\Nw(s)\subset \Delta$ for any moment polytope $\Delta$
described in \S\,2.
Suppose $\gamma\subsetneq\Delta$ is a face of $\Delta$ such that $\gamma\cap\Nw(f) =\varnothing$. 
%One can check that 
The point is to prove that 
any geometry $(M(t),g_1)$ (such as the one in \S\,3) with $t\in \gamma$ is Einstein. So
either $\gamma=\varnothing$, or the set $T$ of Einstein
geometries at infinity contains the whole face $\gamma$ of $\Delta$, that is, $\Delta\ne \Delta_{\min}$.
Therefore, $\Nw(f)=\Delta_{\min}$.
\end{proof}

Now we use the theory of  systems of rational algebraic Laurent equations, developed by A.G.Kushnirenko and D.N.Bernshtein (see e.g., \cite{D.N.Ber}).
(The latter approach via intersections of algebraic cycles on toric varieties is well known.)
It follows from \cite {D.N.Ber} that
$ \varepsilon \le (d-1)! \, V $,
 where $ V $ is the volume of the Newton polytope $\Nw(f)$.
 By Proposition~\ref{PROP:Del_min=Nw}, $ \Nw(f) = \Delta _{\min} $. Hence $ \varepsilon \le \nu $.

% From the results of \cite {D.N.Ber} 
% one can easily derive 
% easily implies that 
 
\smallskip

We now prove the inequality $\nu\le P_{d-1}(3)$, where $P_k$ in $k$-th Legendre polynomial, that is, $P_n(z)=\frac{1}{2^n\,n!}\, \frac{d^n}{dz^n}(z^2-1)^n$.
Using the generating function ${\frac {1}{\sqrt {1-2zw+{w}^{2}}}}=\sum_{k=0}^{\infty} P_k(z)w^k$,
we can write
$$
\begin{gathered}
\sum _{d=1}^{\infty }P_{{d-1}}(3)\,{w}^{d-1} =
{\frac {1}{\sqrt {1-6\,w+{w}^{2}}}}=
\\
1+3\,w+13\,{w}^{2}+63\,{w}^{3}+
321\,{w}^{4}+1683\,{w}^{5}%+8989\,{w}^{6}+48639\,{w}^{7}
%   +265729\,{w}^{8}
%  +1462563\,{w}^{9}  \\  +O\left ({w}^{10}\right ) .
%
+O\left ({w}^{6}\right ). %+O\left ({w}^{8}\right ).
\end{gathered}
$$

The degenerate permutohedron with vertex $p\in \RR^{n+1}$ is the convex hull of the points in Euclidean space obtained from a single point $p$ by all permutations of coordinates.

\begin{LEM}
If $n\ge1$ and $z\ge1$, then
$P_{n}(z)/n!$ is 
the volume of the $n$-dimensional degenerate permutohedron\footnote{The Minkowski sum $\frac{z+1}{2}S+\frac{z-1}{2}S'$ of opposite simplices $S$ and $S'=-S$.} with  vertex 
$(\frac{z+1}{2},0,\dots,0,\frac{1-z}{2})\in\RR^{n+1}$.  
\end{LEM}

\begin{proof}[Proof]
The lemma follows from \cite {Po}, the proof of Theorem 16.3,(8) and Theorem 3.2.
\end{proof}

\smallskip

Let us denote by $ \Pi $ the degenerate permutohedron
with a vertex\footnote{This polytope has $2^d-2=2,6,14,\ldots$ facets, just as the classical (non-degenerate) permutohedra.}  
$$
(2,0,\ldots ,0,-1) \in \RR^d
%(1,0,\ldots ,0,-2) \in \RR^d.
$$
By construction of moment polytopes $\Delta$ (Section~2), we have $\Delta_{\min}\subset \Delta\subset \Pi$, and
the volume of the polytope $ \Pi $ is given by
${V_\Pi
%\opn{vol}(\Pi )
= P_{{d-1}}(3)/(d-1)!} $.
Hence, $\nu\le P_{d-1}(3)$.
% A rough estimate $\nu \le (d-1)! \, V_\Pi % \opn{vol}(\Pi )
% $ follows immediately from the inclusion  $ \Delta _{\min} \subset \Pi $.
Finally,
$$
\varepsilon \le \nu
\le P_{{d-1}}(3)< (3+2\sqrt2)^{d-1} < 6^{d-1}.
$$
The remaining claims of Theorem~\ref {THM:e-nu} follow from \cite[Theorem B]{D.N.Ber}.

%\end{proof}

%(Fig.1).
\begin{figure}
\includegraphics{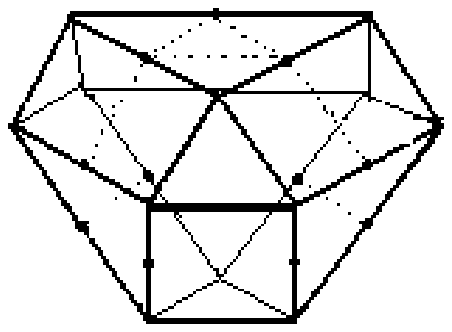}
\caption{For $ d = 4$ the enveloping  permutohedron $ \Pi $ of polytopes $ \Delta $ has %multiplied by $ {(d {-} 1)!} $ 
normalized volume $P_{d-1}(3)=63$ and $2^d-2=14$ facets.}
\end{figure}

\smallskip

The numbers $P_n(3)=1,3,13,63,321,\ldots$ are also known as central Delannoy numbers, that is, $P_n(3)$
%counts the lattice paths running from (0; 0) to (n; n) that use the steps (1; 0), (0; 1), and (1; 1).
counts the number of the paths in $\RR^2$ from $(0, 0)$ to $(n, n)$  that use the steps $(1, 0)$, $(0, 1)$, and $(1, 1)$.
These numbers form the diagonal of the symmetric array $(d_{m,n})$ introduced by H. Delannoy (1895) in the same way, so that
$
{d_{m,n}=d_{m-1,n}+d_{m,n-1}+d_{m-1,n-1}}.
$ 
See, e.g., \cite[6.3.8]{Stanley2}.  

\begin{COR*} 

We have the upper bound of the $(d-1)$-th central Delannoy number  for the  normalized volume $\nu$ of the moment polytope.
Here $d>1$ is the number of the irreducible submodules in the isotropy $H$-module $\frak g/\frak h$.
 
\end{COR*}

%\bibitem{Stanley2} R. P. Stanley, Enumerative Combinatorics, Vol. 2, Cambridge University Press, 1999.

\begin{EXAMPLES}{} The inequality $ \varepsilon < \nu$ holds for any space $G/H$ with $T\ne \varnothing$, e.g., for the spaces in Examples~\ref{EXAM:31}, ~\ref{EXAM:32}.
 Consider $SU(p)$ for small $p$ as a homogeneous space $G/H$ in Example~\ref{EXAM:31}.
Let $p\in \{2,3\}$. Then $T$ is a finite set.
One can check (by an examination of the
solutions at infinity only) that $\varepsilon = \nu- |T|$.
The Newton polytopes are just the same as for %the spaces
$SU(3)/T^2$ and $E_8/(A_2)^4$ (Examples~\ref{EXAM:22}, \ref{EXAM:23}), and $\nu \in \{4,23\}$.
Hence, $\varepsilon = \nu-p-1= 1$ and $19$ respectively.
\end{EXAMPLES}

\begin{EXAMP}{} 
Let $G/H $ be the $196$-dimensional K\"ahler homogeneous
space $E_8/T^1\cdot A_1\cdot A_6$. % called Louis XIV. (\cite{14}).
Then $d=4$ and $\varepsilon = \nu = 20$.
This is less than $1/3\,P_3(3)=21$.
\end{EXAMP}

\bigskip



\section*{Appendix. Case of K\"ahler homogeneous space with $b_2=1$}

\let\phim m

Consider now a K\"ahler homogeneous space $G/H$ with the second Betti number $b_2=1$.
Assume that the isotropy $H$-module $\frak g/\frak h$ is split into $ d>1 $ irreducible submodules.

\begin{LEM}%\label{LEM:} 

Given a K\"ahler homogeneous space $G/H$ with $d>1$ of a simple Lie group $G$, then $2^{-b_2(G/H)}\nu(G/H) \in \ZZ$.

\end{LEM}

\begin{proof}[Idea of proof] $2^{b_2(G/H)} = [\ZZ^d:L]$,
where $L \subset \ZZ^d$ is the subgroup generated by vertices of the polytope $\Delta$. 
(Remark that $\Delta = \Delta_{\min}=\Delta_{\max}$).
\end{proof}

Let, moreover, $b_2(G/H)=1$.
Then $2\le d \le 6$. %$d \in \{2,\dots,6\}$.  
For $d=2$ the polytope $\Delta$ is the segment with ends $e_2$ and $2e_1-e_2$. 
For $3\le d \le 6$ the vertices of the polytope $\Delta $ are the points 
$e_i+e_j-e_k \in \RR^d$
with $1\le i,j,k \le d$, $i\ne k$,  
$j\ne k$, $i\pm j\pm k =0$. Here $e_1=(1,0,\ldots,0),\dots,e_d=(0,\dots,0,1)$.
Using MAPLE, one can triangulate these polytopes and find their normalized volumes $\nu=\nu(G/H)$.
Thus, we obtain:
\begin{CLAIM*} 
$
2^{-1}\nu \in\{1,3,10,41,172\}.
$
\end{CLAIM*}

The following table gives  some information
about polytopes $\Delta$ corresponding to K\"ahler homogeneous spaces $G/H$ with the second Betti number $b_2=1$ and $d>1$.

For completeness, we find the volume of a similar $(d-1)$-dimensional polytope with $d=7$. 
Here $f$ is the number of facets of $\Delta$, and $\phim$ the number of all  faces $\gamma$ of $\Delta$ with $0<\dim(\gamma)<d-1$ 
(which we call \textbf{marked}) 
that NOT satisfy conditions of Test~1 or Test~2 of \cite[\S 7.1]{2007}.
A marked face $\gamma$ is not a vertex.
\footnote{
Let $G/H$ be a K\"ahler homogeneous space of a simple Lie group $G$, and let $d>2$. Then the set of vertices of $\Delta$ is $\{e_i+e_j-e_k : 1\le i,j,k\le d, i\ne k,j\ne k, [i,j,k]\ne0\}$. In this case, conditions of Tests~1 and~2 for a $k$-dimensional face $\gamma$, $0<k<d-1$ can be simplified as follows: \\ 1)~$\gamma $ is a pyramid with apex $a$ and base $B$ such that if $e_i\in \gamma$, then either $e_i=a$ or $e_i\in B$; \\ 2)~$\gamma$ is a `$k$-dimensional octahedron' with vertices $e_{i_0} + v_p, e_{i_0}-v_p$, $p=1,\dots,k$, for some linearly independent vectors $v_p\in\RR^d$; the face $\gamma$ contains no points $e_i$ with $i\ne i_0$ (then $\gamma$ is the intersection of all faces $\beta\ni e_{i_0}$ of $\Delta$). For $b_2(G/H)=1$, $d>2$ there are $[d/2]$ faces $\gamma$ satisfying~2).   
}
Moreover, one can prove that $\gamma$ is not an edge, so $1<\dim(\gamma)<\dim(\Delta)$. 
$$
\begin{array}{c|llllll|llll}
d   & & 2 & 3 & 4  & 5  & 6  & 7 \\
f   & & 2 & 4 & 7  & 16 & 36 & 100 \\
\nu & & 2 & 6 & 20 & 82 & 344 & 1598 \\ \hline \varepsilon
    & &\nu&\nu&\nu & 81 &  ?  & - \\ \delta
	& & 0 & 0 & 0  &  1 &  ?  & -  \\ \phim
	& & 0 & 0 & 3  & 13 &  40 & 
\end{array}
$$
We write also the known numbers $\varepsilon$ and $\delta = \nu-\varepsilon$. 
By \cite{2007}, if $m=0$, then $\nu=\varepsilon$.
The first non-trivial case is $d=4$.

\medskip

For $d\le 5$ all the \textit{positive} solutions of the  algebraic Einstein equations are known.
% We write also the known numbers $\varepsilon$ and $\delta = \nu-\varepsilon$. 
%In particular, the case $d=5$ is recently studied 
In the case $d=5$ they calculated by I.Chrysikos and Y.Sakane \cite{SACOS}.
% \footnote{Ioannis Chrysikos, Yusuke Sakane. On the classification of homogeneous Einstein metrics on generalized flag manifolds with $b_2(M)=1$. arXiv:1206.1306v1}
They also prove that all the complex solution are isolated ($d=5$).

\medskip

{\bf Remark (the case $d=4$, $\nu=20$). }
Here $\Delta$ is a three-dimensional polytope with three marked faces $\gamma,$ namely, a trapezoid $\gamma_1$, a parallelogram $\gamma_2$, and a pentagon $\gamma_3$. To prove that $\nu=\varepsilon$ one can associate with each marked face $\gamma$ a complex hypersurface $s_{\gamma}(x_1,\dots,x_4)=0$ in $(\CC \setminus 0)^4$ and check that it is non-singular.  
Here $s(x)$ is the above Laurent polynomial (scalar curvature), and $s_{\gamma}(x)$ is the sum of all monomials of  $s(x)$ whose vector exponents belong to $\gamma$. See \cite[\S1.7.2]{2006,2007}. This is essentially a two-dimensional problem ($2=\dim(\gamma_i)$), i.e., we must check that a plane curve is non-singular. It is easy to prove this for $\gamma=\gamma_1$, but for $\gamma_k$, $k=2,3$ the problem reduces to 
$D_k[s]\ne0$, where $D_k[s]$ is a homogeneous polynomial 
 (a $k$-monomial)
in coefficients of $s(x)$ with $\deg  (D_k) = k$. The coefficients of $s(x)$ depend on $G/H$. There are four K\"ahler homogeneous spaces $G/H$ with $b_2=1$ and $d=4$, namely, the spaces
$
{E_8/T^1\cdot A_1\cdot A_6},\quad
{E_8/T^1\cdot A_2\cdot D_5},\quad
{E_7/T^1\cdot A_1\cdot A_2\cdot A_3},\quad
{F_4/T^1\cdot \widetilde A_1\cdot A_2}
$
(we use the Dynkin's notation $\widetilde A_1$ for the three-dimensional subgroup associated with a short simple root).
The corresponding scalar curvature polynomials $s(x)$ are computed by A.Arvanitoyeorgos and I.Chrysikos (arXiv:0904.1690). For each of them one can check that $D_2[s]D_3[s]\ne0$. This proves that $\varepsilon=\nu$ for $d=4$.

\begin{figure}[h]
\includegraphics[scale=0.7]{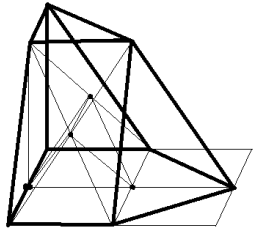}
\caption{The $3$-dimensional polytope $\Delta$ with $7$ facets corresponding to four K\"ahler homogeneous spaces with $b_2=1$ and $d=4$}
\end{figure}

\medskip

{\bf The case $d=5$, $\nu=82$.}
There is a unique K\"ahler homogeneous space $G/H$ with $b_2(G/H)=1$
and $d=5$, namely, the space 
${G/H=E_8/T^1\cdot A_3\cdot A_4}.$
By \cite[\S3, the text after eq.(25)]{SACOS} it implies that 
the algebraic Einstein equation has, up to scale,  $81$ complex solutions, corresponding to roots of some polynomial $(x_5-5)h_1(x_5)$ of degree $81$ in one variable $x_5$.
There exist $6$ positive solutions \cite[Theorem A]{SACOS};
in particular, the root $x_5=5$ corresponds to a unique, up to scale, invariant K\"ahler metric on $G/H$.
Using MAPLE, one can check that the polynomial $h_1(x_5)$ has $80$ simple roots, and $81$ solutions of Einstein equation are distinct
(moreover, it has $30$ real roots, and Einstein equation has $31$ real solutions).
Thus
$$
\nu - \varepsilon = 82-81 = 1.
$$

We prove %directly
independently
that $\nu > \varepsilon $.
Let $s(x_1, \dots ,x_5)$ be the scalar curvature
of a invariant metric $g_x$, as above. It is a Laurent polynomial in $x_i^{-1}$.
We claim that there exists a limit
$$
s_{\infty}(x_1, \dots ,x_5) =
\lim_{t\,\to\,+ \infty } s(t^{2}\,{x_{1}}, t^{4}\,{x_{2}},
t^{3}\,{x_{3}}, t\,{x_{4}}, t\,{x_{5}}),
$$
and the homogeneous function $s_{\infty}$ depends essentially on $3$ variables.
Indeed, it follows from Proposition~\ref{PROP:Del_min=Nw} and the above description of the Newton polytope
$\Delta $ that
$$
\begin{aligned}
s_{\infty} 
=&
- \frac{\scriptstyle [1,4,5]}{2}{\frac {{x_{1}}}{{x_{4}}\,{x_{5}}}}
- \frac{\scriptstyle [1,3.4]}{2}{\frac {{x_{3}}}{{x_{1}}\,{x_{4}}}}
\\ &
- \frac{\scriptstyle [2,3,5]}{2}{\frac {{x_{2}}}{{x_{3}}\,{x_{5}}}}
- \frac{\scriptstyle [1,1,2]}{4}{\frac {{x_{2}}}{{x_{1}}^{2}}}.
\end{aligned}
$$
Then there is a two-dimensional face of the polytope $\Delta$, the  parallelogram $P$   with vertices $e_4+e_5-e_1$, $e_1+e_4-e_3$, $2e_1-e_4$, and $e_3+e_5-e_2$. The face $P$ is  orthogonal to the vector
  $$
\mathbf f = (2,4,3,1,1) . %\in \nabla
  $$
\begin{comment}
Then there the parallelogram $P$ with vertices $e_4+e_5-e_1$, $e_1+e_4-e_3$,
$2e_1-e_4$, and $e_3+e_5-e_2$ is a two-dimensional face of the polytope $\Delta$,
orthogonal to the vector
$$
\mathbf f=(2,4,3,1,1).
$$
\end{comment}
According to \cite[Proposition 7]{SACOS} we have
$
[1, 1, 2]=12,
[1, 2, 3]=8,
[1, 3, 4]=4,
[1, 4, 5]=4/3,
[2, 2, 4]=4,
[2, 3, 5]=2
$.
Then
\begin{comment} *b**************************************************
$
D_P = \det \left( \begin{smallmatrix}
[1,4,5]/2 &\,\,
[1,3,4]/2 \\[1ex]
[2,3,5]/2 &\,\,
[1,1,2]/4
\end{smallmatrix}\right)
= \det \left(\begin{smallmatrix} 4/6& 2\\[1ex] 1 & 12/4  \end{smallmatrix}\right)
=0,
$
%%
\end{comment}
%%
%%%*************************************************************e%%%
the product of monomials, corresponding to each pair of opposite vertices of $P$,
coincides with $2 {\frac {{x_{2}}}{{x_{1}}\,{x_{4}}\,{x_{5}}}}$, and $s_{\infty }$
can be represented as
$$
s_{\infty} = z_0(1+z_1+z_2+z_1z_2) = z_0(z_1+1)(z_2+1),
$$
where
$z_0=-\frac {{x_{2}}}{{x_{3}}\,{x_{5}}}$.
Since for $z_1=z_2=-1$ we have
$$
s_{\infty } = ds_{\infty } =0,
$$
the complex hypersurface
$s_{\infty }(x)=0$ has a singular point $x$ with $\prod x_i\ne0$.
% By \cite[Proposition ..]{2006,2007}
By \cite[\S1.7.2]{2006,2007} 
this implies that $\nu- \varepsilon >0$.

Note that $\Delta$ has $\phim-1=12$ marked faces, other than $P$;
namely, $6$ three-dimensional faces with normal vectors
$$
\mathbf f =
%npyram_faces3 :=
(1, 2, 3, 4, 5), (1, 2, 1, 2, 1), (2, 1, 1, 2, 0), %(1, 1, 0, 1, 1), removed!
 (1, 0, 1, 0, 1),(1, 2, 2, 1, 0), (1, 2, 1, 0, 1),
%(1, 2, 2, 1, 0), (1, 0, 1, 0, 1), (1, 2, 1, 0, 1), transposition!
$$
and $6$ parallelograms defined by the following normal vectors
(such as $\mathbf f=(2, 4, 3, 1, 1)$):
$$
\mathbf f =
(1, 1, 2, 2, 3), (1, 2, 3, 4, 4), (1, 2, 2, 3, 4), (2, 4, 5, 3, 1),
(5, 3, 2, 6, 1), %,
(3, 1, 2, 2, 1)
%, (1, 1, 1, 0, 1), (1, 1, 1, 1, 0)
$$
The corresponding $12$ complex hypersurfaces are non-singular.

\bigskip

%\NULLEMPTY{\input ./THAT50/Localization.tex}

{\bf Additional remarks (the case $d=5$). }
Consider $(z_0,z_1,z_2)=(1,-1,-1)$ as a point $p$ in
the four-dimensional toric variety $\Delta^{\CC}$. Let $O\subset\Delta^{\CC}$  be the orbit of the group $(\CC \setminus 0)^5/\CC^\times$ trough $p$. The closure of $O$ is the two-dimensional toric subvariety $P^{\CC}$.
The point $p\in O$ is a solution at infinity (in the sence of \S7)  of the algebraic Einstein equation.

\smallskip

Our example is excellent as the following lemma show.

\begin{LEM}

We claim now that $\Delta^{\CC}$ is smooth at each point $q\in O$. 
Moreover, assuming
$\phi:\Delta^{\CC}\to\mathbb P^{N-1}(\mathbb C)$ 
be the natural map into the complex projective space $\mathbb P^{N-1}(\mathbb C)$, $N=\#(\ZZ^5\cap\Delta) $, then $\phi^{-1}(\phi(q))=\{q\}$, and $\phi(\Delta^{\CC})$ is smooth at %the point $\phi(q)$. 

\end{LEM}

We will apply the localization along $O$ to prove
that the point $p$ is an isolated solution (at infinity) 
with the multiplicity $1$ of the algebraic Einstein equation.

\begin{proof}[Proof]
Let $v_0,v_1,v_{12},v_2$
are vertices of the parallelogram $P$, so $v_0+v_{12}=v_1+v_2$, and 
$$
\begin{aligned}
u_1 & = -e_1+e_4+e_5 & = v_1, &
\\ 
u_2 & = -e_2 + 2e_4+2e_5 & = 2v_1+v_2, &
\\ 
u_3 & = -e_3+2e_4+e_5 & = v_1+v_{12}, &
\quad\,
u_4=e_4,\quad\,
u_5=e_5.
\end{aligned}
$$
The set of vectors $\{u_i : i=1,\dots,5\}$, 
and, hence, $\{v_0,v_1,v_2,u_4,u_5\}$ are basises in $\ZZ^5=\bigoplus \ZZ e_i$.
Let $\pi(a):=(a_4,a_5)$ for each $a=\sum a_i u_i \in \ZZ^5$.
We prove, that $\pi(\ZZ^5\cap \Delta)$ generates the semigroup $\ZZ_+^2$.
The face $P$ of $\Delta$ is the intersection of two facets with normal vectors $\mathbf f_i$, $i=1,2$, so that
$$
\mathbf f = (2,4,3,1,1) = (1, 2, 2, 1, 0) + (1, 2, 1, 0, 1) =\mathbf f_1 + \mathbf f_2 .
$$
For any $a\in \ZZ^5\cap \Delta$ we have 
$a_4 = \langle \mathbf f_1,a\rangle\ge0$, and
$a_5 = \langle \mathbf f_2,a\rangle\ge 0$.
Then $\pi(a)\in \ZZ_+^2$. This prove the assertion, since $\pi(e_4) = (1,0)$, $\pi(e_5)=(0,1)$, $e_4,e_5\in\Delta$. 
The lemma follows
\end{proof}

Now let $(z_0,z_1,z_2,y_1,y_2)$ be coordinates on $(\CC \setminus 0)^3\times \CC^2$. Assume that for $y_1y_2\ne0$ 
\[
- {\displaystyle \frac {{x_{2}}}{{x_{3}}\,{x_{5}}}} ={z_{0}}, 
\, - {\displaystyle \frac {2}{3}} \,{\displaystyle \frac {{x_{1}}
}{{x_{4}}\,{x_{5}}}} ={z_{0}}\,{z_{1}}, \, - 3\,{\displaystyle 
\frac {{x_{2}}}{{x_{1}}^{2}}} ={z_{0}}\,{z_{2}}, \, - 2\,
{\displaystyle \frac {{x_{3}}}{{x_{1}}\,{x_{4}}}} ={z_{0}}\,{z_{1
}}\,{z_{2}}, \,{\displaystyle \frac {1}{{x_{4}}}} ={y_{1}}, \,
{\displaystyle \frac {1}{{x_{5}}}} ={y_{2}},
\]
so %that for $y_1y_2\ne0$ we have
\[
 \left\{  \! {x_{3}}={\displaystyle \frac {3}{4}} \,
{\displaystyle \frac {{z_{0}}^{2}\,{z_{1}}^{2}\,{z_{2}}}{{y_{1}}
^{2}\,{y_{2}}}} , \,{x_{2}}= - {\displaystyle \frac {3}{4}} \,
{\displaystyle \frac {{z_{0}}^{3}\,{z_{1}}^{2}\,{z_{2}}}{{y_{1}}
^{2}\,{y_{2}}^{2}}} , \,{x_{4}}={\displaystyle \frac {1}{{y_{1}}}
} , \,{x_{1}}= - {\displaystyle \frac {3}{2}} \,{\displaystyle 
\frac {{z_{0}}\,{z_{1}}}{{y_{1}}\,{y_{2}}}} , \,{x_{5}}=
{\displaystyle \frac {1}{{y_{2}}}}  \!  \right\}. 
\]
Then 
\begin{multline*}
s= \\
 - {\displaystyle \frac {16}{9}} \,{\displaystyle \frac {{y_{1}}
^{3}\,{y_{2}}^{4}}{{z_{0}}^{6}\,{z_{1}}^{4}\,{z_{2}}^{2}}}  + 
{\displaystyle \frac {{\displaystyle \frac {16}{9}} \,{y_{1}}^{4}
\,{y_{2}}^{2}}{{z_{0}}^{5}\,{z_{1}}^{4}\,{z_{2}}^{2}}}  - 
{\displaystyle \frac {32}{3}} \,{\displaystyle \frac {{y_{1}}^{3}
\,{y_{2}}^{2}}{{z_{0}}^{4}\,{z_{1}}^{3}\,{z_{2}}^{2}}}  + 
{\displaystyle \frac {{\displaystyle \frac {16}{9}} \,{y_{1}}^{2}
\,{y_{2}}^{2}}{{z_{0}}^{3}\,{z_{1}}^{3}\,{z_{2}}}}  - 
{\displaystyle \frac {32\,{y_{1}}^{2}\,{y_{2}}^{2}}{{z_{0}}^{3}\,
{z_{1}}^{2}\,{z_{2}}}}  + {\displaystyle \frac {{\displaystyle 
\frac {80}{3}} \,{y_{1}}^{2}\,{y_{2}}}{{z_{0}}^{2}\,{z_{1}}^{2}\,
{z_{2}}}}  - {\displaystyle \frac {8}{3}} \,{\displaystyle 
\frac {{y_{1}}\,{y_{2}}^{2}}{{z_{0}}^{2}\,{z_{1}}}}  \\
\mbox{} + {\displaystyle \frac {4\,{y_{1}}^{2}}{{z_{0}}\,{z_{1}}
\,{z_{2}}}}  + {\displaystyle \frac {{\displaystyle \frac {4}{9}
} \,{y_{1}}^{2}}{{z_{0}}\,{z_{1}}}}  - {\displaystyle \frac {80}{
3}} \,{\displaystyle \frac {{y_{1}}\,{y_{2}}}{{z_{0}}\,{z_{1}}}} 
 + {\displaystyle \frac {{y_{2}}^{2}}{{z_{0}}}}  + 
{\displaystyle \frac {{\displaystyle \frac {4}{9}} \,{y_{2}}^{2}
}{{z_{0}}\,{z_{1}}}}  + 8\,{y_{1}} - {\displaystyle \frac {8}{3}
} \,{\displaystyle \frac {{y_{1}}}{{z_{1}}}}  + 4\,{y_{2}} + {z_{
0}}\,{z_{1}} + {z_{0}}\,{z_{1}}\,{z_{2}} \\
\mbox{} + {z_{0}} + {z_{0}}\,{z_{2}}, 
\end{multline*}
$s$ is a polynomial in $y_1$, $y_2$, and 
$$
s= z_0 + z_0z_1+z_0z_2 +z_0z_1z_2 + 8y_1 - \frac83\frac{y_1}{z_1} +4y_2  + [2],
$$
where $[2]$ denotes the terms with degree $\ge 2$.
Similarly, for $s_i = x_i\partial s/\partial x_i$, $i=1,\dots,5$ we have
$$
\begin{aligned}
&
s_1= z_0z_1 -2z_0z_2 -z_0z_1z_2 + \frac83\frac{y_1}{z_1} + [2],
&
s_2= + z_0 + z_0z_2 - \frac83\frac{y_1}{z_1} + [2],
\\ &
s_3= -z_0 + z_0z_1z_2 + \frac83\frac{y_1}{z_1} + [2],
&
s_4= -z_0z_1 - z_0z_1z_2 - 8y_4 + [2],
\\ &
s_5= -z_0 - z_0z_1 - 4y_5 + [2],
\end{aligned}
$$
where $[2]$ denotes $(y_1^2, y_1y_2, y_2^2)$.
Computing the matrix $J= \frac{\partial(s_1,s_2,s_3,s_4,s_5)}{\partial(z_1,z_2,y_1,y_2)}$, setting $z_0=1,z_1=z_2=-1$, $y_1=y_2=0$, adding the row $(d_1,\dots,d_5)$ of dimensions $d_i=\dim(\frak m_i)$, and finding the determinant, we obtain
$$
\left|\begin{array}{rrrrr}
d_1&d_2&d_3&d_4&d_5\\ 2&0&-1&0&-1\\ -1&1&-1&1&0 \\ -8/3&8/3&-8/3&-8&0 \\ 0&0&0&0& -4
\end{array}\right| %1&-1&1&3
=\frac{128}{3} (d_1+3d_2+2d_3) >0.
$$ 
Then the solution $p\in\Delta^{\CC}$ of the algebraic Einstein equation with local coordinates
$z_0=1,z_1=z_2=-1, y_1=y_2=0$ is isolated, and non-degenerate.

\bigskip

 {\bf The case $d=6$, $\nu = 344$.}
There is a unique K\"ahler homogeneous space $G/H$ with $b_2(G/H)=1$
and $d=6:$ 
$$
G/H = E_8/T^1\cdot A_4\cdot A_2\cdot A_1.
$$
The corresponding $5$-dimensional polytope $\Delta $ in $\RR^6$
has $36$ facets, i.e., $4$-dimensional faces. Each
of them can be defined by the orthogonal vector $\mathbf f=(y_1, \dots ,y_6)$
such that 
$\gcd(y_1, \dots ,y_6)=1$, and
$y_i\ge0$; then $\langle \mathbf f,x \rangle\,\ge 0$ for any $x \in  \Delta $.

For example, the vector $\mathbf f=(1, 2, 3, 4, 5, 6)$ is orthogonal to the facet
with $9$ vertices
$$
\def\1#1,#2,#3{e^{#3}_{#1 #2}}
\begin{array}{llllllll}
\11,1,2,&\11,2,3,&\11,3,4,&\11,4,5,&\11,5,6,\\[1ex]
        &\12,2,4,&\12,3,5,&\12,4,6,& \\[1ex]
        &        &\13,3,6,& &
\end{array}
$$
where $e^k_{ij}=e_i+e_j-e_k$. We write all the facets:

1) \enskip  $16$ four-dimensional simplices with normal vectors
{\footnotesize
\begin{verbatim}
 [1, 2, 2, 1, 2, 1], [1, 2, 3, 2, 3, 2], [1, 2, 1, 1, 2, 1], [1, 1, 1, 2, 2, 1],
 [2, 2, 1, 1, 1, 1], [3, 2, 3, 4, 1, 2], [1, 1, 1, 2, 1, 2], [1, 1, 2, 2, 1, 1],
 [1, 1, 2, 2, 1, 2], [2, 1, 1, 1, 1, 2], [3, 2, 5, 4, 3, 6], [1, 1, 2, 1, 1, 2],
 [2, 2, 3, 1, 1, 3], [5, 2, 3, 4, 5, 6], [2, 1, 1, 1, 2, 2], [5, 4, 3, 2, 7, 6];
\end{verbatim}
}

2) \enskip $8$ four-dimensional  pyramids with normal vectors
{\footnotesize
\begin{verbatim}
 [3, 4, 1, 2, 5, 2], [1, 2, 3, 4, 5, 4], [5, 2, 3, 4, 1, 6], [1, 2, 3, 2, 1, 2],
 [3, 2, 1, 2, 1, 2], [1, 2, 3, 4, 3, 4], [3, 2, 1, 4, 3, 2], [1, 2, 3, 2, 3, 4];
\end{verbatim}
}

3) \enskip  $5$ other facets  with normal vectors with positive entries:
{\footnotesize
\begin{verbatim}
 [1, 2, 3, 4, 5, 6], [3, 2, 1, 4, 1, 2], [1, 2, 3, 4, 3, 2], [1, 2, 1, 2, 3, 2],
 [1, 2, 1, 2, 1, 2];
\end{verbatim}
}

4) \enskip  $7$ facets  with normal vectors with non-negative entries:
{\footnotesize
\begin{verbatim}
 [1, 2, 1, 0, 1, 2], [2, 1, 1, 2, 0, 2], [1, 2, 2, 1, 0, 1], [1, 1, 0, 1, 1, 0],
 [1, 0, 1, 0, 1, 0], [1, 2, 1, 2, 1, 0], [1, 2, 3, 2, 1, 0];
\end{verbatim}
}

Facets 1) and 2) are not marked faces.
E.g., simplices 1) and its sub-faces satisfy \cite[Test 7.1]{2007}.
% The pyramids 2) may be replaced by its maximal non-pyramidal faces.

Facets 3) and 4) are marked faces. There are $13$ three-dimensional
and $15$ two-dimensional marked faces.

We get $13$ vectors, orthogonal to three-dimensional marked faces (each of them is proportional to the sum of two distinct vectors 2)-4)): 
% Each of them, up to scale, can be represented as a
 % sum over the above $20=8+5+7$
% Each of them is proportional to a linear combination of vectors 2)-4) with coefficients $0$ and $1 :$ 
% Each of them is proportional to a linear combination of vectors 2)-4) with coefficients $0$ and $1 :$
{\footnotesize
% mk3:=
\begin{verbatim}
 [1, 2, 1, 2, 1, 1], [1, 2, 1, 1, 1, 2], [5, 3, 2, 6, 1, 4], 
 [4, 3, 1, 5, 2, 2], [7, 3, 4, 6, 1, 8], [4, 5, 1, 3, 6, 2], 
 [1, 2, 3, 4, 5, 5], [2, 4, 5, 3, 1, 1], [2, 4, 3, 1, 1, 3], [1, 2, 2, 2, 1, 0], 
 [1, 1, 2, 1, 1, 0], [2, 1, 1, 1, 2, 0], [2, 3, 1, 3, 2, 0]
\end{verbatim}
}

In the two-dimensional case we obtain:

(a) \enskip $6$ parallelograms with normal vectors
{\footnotesize
\begin{verbatim}
 [5, 5, 2, 7, 3, 2],    [3, 6, 8, 5, 2, 3], [5, 5, 2, 3, 7, 2],
 [3, 6, 7, 10, 13, 12], [3, 4, 6, 3, 2, 1], [3, 6, 6, 5, 2, 1]
\end{verbatim}
}

(b) \enskip $9$ parallelograms with normal vectors
{\footnotesize
\begin{verbatim}
 [3, 6, 7, 8, 5, 2],  [8, 3, 5, 6, 2, 8], [5, 7, 2, 3, 7, 4],
 [3, 6, 7, 6, 9, 12], [7, 5, 2, 9, 5, 4], [5, 7, 2, 5, 9, 4],
 [3, 4, 7, 8, 11, 10],[3, 2, 1, 3, 1, 2], [1, 2, 3, 4, 4, 4]
\end{verbatim}
}%
Each of parallelograms listed in (a) and (b)
(with the exception of two last entries in (b))
belongs exactly to $3$ facets.

We claim now, that $6$ marked parallelograms (a) corresponds to singular complex hypersurfaces
as above (consequently $\varepsilon<\nu$), and $9$ marked parallelograms (b) corresponds to non-singular hypersurfaces.
For the proof, one can calculate $6+9$ determinants $\left|\begin{smallmatrix}a&b\\ b'&a'\end{smallmatrix}\right|$,
where $a,a',b,b'$ are some coefficients of $s(x_1,\dots,x_6)$,
using equalities
(\cite[Prop. 13]{SACOS}):
$
[1, 1, 2] = 8, [1, 2, 3] = 6, [1, 3, 4] = 4, [1, 4, 5] = 2, [1, 5, 6] = 1, [2, 2, 4] = 6, [2, 3, 5] = 2, [2, 4, 6] = 2, [3, 3, 6] = 2
$.

Thus we may unmark $9$ of $40$ marked faces.

\begin{COR*}
The hypothesis that all complex solutions of the algebraic Einstein equation on
$
G/H = E_8/T^1\cdot A_4\cdot A_2\cdot A_1
$
are isolated reduces to examination of $31 = 40-9$ cases, corresponding to $12$ four-dimensional, 
$13$ three-dimensional, and $6$ two-dimensional faces of the polytope $\Delta $.
\end{COR*}

In each case we may unmark the $k$-dimensional face, if the corresponding complex hypersurface is non-singular (this is the $k$-dimensional problem, $k<5$); otherwise we must examine Einstein equation in a neighborood $U \subset \Delta^{\CC}$ of the 'solution at infinity' defined by each singular point (cf. \S\ref{sect:7}).



\medskip

I am grateful to all participants of the Postnikov seminar of Moscow State University, listening to the presentation of this work May 20, 2009 
and March 23, 2011.

%\hfil <<█═╓╔НАЛ, ╒═╛ ╞╝╜ОБ╜╝, ё╝А╞╝╓═?>> (▒╝╜ ▐╝╞╝╒═)


\begin{thebibliography}{99} 
\label{bibliography}

\bibitem{Al-Ki} D. V. Alekseevsky and B. N. Kimel'fel'd, Structure of homogeneous Riemann spaces with zero Ricci
curvature, Func. Anal. Appl. 9 (1975), 97-102.

\bibitem{AB}   A.L.Besse. Einstein Manifolds, Springer-Verlag, Berlin, Heidelberg, New York, 1987.

\bibitem{D.N.Ber} D. N. Bernshtein, The number of roots of a system of equations, Funktsional. Anal. i
Prilozhen. 9:3 (1975), 1-4; English transl., Funct. Anal. Appl. 9:3 (1975), 183-185.


\bibitem{BWZ}  C.Boehm-M.Wang-W.Ziller, A variational approach for compact homogeneous Einstein manifolds, GAFA 14 (2004), 681-733

\bibitem{B-G}  Charles P. Boyer and Krzysztof Galicki, Sasakian geometry, Oxford Mathematical Monographs, Oxford University Press, Oxford, 2008.


\bibitem{Fulton} W. Fulton. Introduction to toric varieties, Princeton University Press, 1993.


\bibitem{2006} M.M.Graev. On the number of invariant Einstein metrics on a compact homogeneous space, Newton polytopes and contractions of Lie algebras.
%International Journal of Geometric Methods in Modern Physics.
Int. J. G. M. Mod. Phys.  Vol. 3, Nos. 5 \& 6 (2006) 1047-1075

\bibitem{2007} 
M. M. Graev. The number of invariant Einstein metrics on a homogeneous space, Newton polytopes and contractions of Lie algebras. Izvestiya RAN: Ser. Mat. 71:2 (2007), 29-88.
English transl., Izvestiya: Mathematics 71:2 (2007) 247-306  

\bibitem{He} Jens Heber .  Noncompact homogeneous Einstein spaces.   Invent. math. 133, 279-352 (1998)  

\bibitem{Jen-2} Gary Jensen, The Scalar Curvature of Left-Invariant Riemannian Metrics, Indiana Univ. Math. J. 20 No. 12 (1971), 1125-1144

\bibitem{Lauret} Jorge Lauret.  Convergence of homogeneous manifolds. arXiv:1105.2082

\bibitem{Po} Alexander Postnikov: Permutohedra, associahedra, and beyond, arXiv:math.CO/0507163.
 
  
\bibitem{Stanley2} R. P. Stanley, Enumerative Combinatorics, Vol. 2, Cambridge University Press, 1999.
 
\bibitem{SACOS}
Ioannis Chrysikos, Yusuke Sakane,
On the classification of homogeneous Einstein metrics on generalized flag manifolds with $b_2(M)=1$.
 arXiv:1206.1306v1 
  
\end{thebibliography}
\end{document}